\documentclass[11pt]{amsart}

\usepackage{amssymb} %
\usepackage{bbm} %
\usepackage{geometry}
\usepackage{graphicx} %
\usepackage{enumitem}
\usepackage[only,llbracket,rrbracket]{stmaryrd}
\usepackage[all]{xy}

\usepackage{hyperref}
\usepackage{xcolor}
\definecolor{darkgreen}{rgb}{0,0.5,0}
\hypersetup{colorlinks,urlcolor=blue,citecolor=darkgreen,linkcolor=blue,linktocpage=false}

\usepackage[abbrev]{amsrefs}

\numberwithin{subsection}{section}
\numberwithin{equation}{section}
\numberwithin{figure}{section}

\theoremstyle{plain}
\newtheorem{theorem}[equation]{Theorem}
\newtheorem{proposition}[equation]{Proposition}
\newtheorem{lemma}[equation]{Lemma}
\newtheorem{corollary}[equation]{Corollary}
\newtheorem{conjecture}[equation]{Conjecture}

\theoremstyle{definition}

\newtheorem{definition}[equation]{Definition}
\newtheorem{example}[equation]{Example}
\newtheorem{notation}[equation]{Notation}

\theoremstyle{remark}
\newtheorem{remark}[equation]{Remark}
\newtheorem{warning}[equation]{Warning}

\newcommand{\A}{\mathbb{A}} %

\newcommand{\F}{\mathbb{F}}
\newcommand{\G}{\mathbb{G}}
\newcommand{\kk}{\Bbbk}
\newcommand{\KK}{\mathbb{K}}

\newcommand{\PP}{\mathbb{P}}
\newcommand{\Q}{\mathbb{Q}}

\newcommand{\Z}{\mathbb{Z}}

\newcommand{\unit}{\mathbbm{1}}
\newcommand{\al}{\alpha}
\newcommand{\be}{\beta}

\newcommand{\ep}{\epsilon}
\newcommand{\ga}{\gamma}

\newcommand{\la}{\lambda}

\newcommand{\om}{\omega}
\newcommand{\phy}{\varphi}
\newcommand{\Si}{\Sigma}
\newcommand{\si}{\sigma}
\newcommand{\te}{\theta}

\newcommand{\inj}{\hookrightarrow}
\newcommand{\ral}{\xrightarrow} %
\newcommand{\rla}{\rightleftarrows}

\newcommand{\bigop}{\bigoplus}

\newcommand{\dual}{\vee}

\newcommand{\op}{\oplus}
\newcommand{\ot}{\otimes}

\newcommand{\sm}{\wedge} %
\newcommand{\x}{\times}

\newcommand{\dfn}{:=}
\newcommand{\hatt}{\widehat}
\newcommand{\tild}{\widetilde}
\newcommand{\ul}{\underline}

\newcommand{\Ab}{\mathrm{Ab}}
\newcommand{\base}{\mathcal{B}} %
\newcommand{\Cat}{\mathbf{Cat}}
\newcommand{\cat}[1]{\mathcal{#1}} %
\newcommand{\Ch}{\mathrm{Ch}}

\newcommand{\hMod}[1]{D(#1)}

\newcommand{\abs}[1]{\lvert #1 \rvert}
\newcommand{\Brack}[1]{\llbracket #1 \rrbracket}

\newcommand{\Def}[1]{\textbf{\boldmath{#1}}} %

\newcommand{\coact}{\mathrm{coact}}
\newcommand{\Cof}{\mathrm{Cof}}
\newcommand{\Conn}{\pi_0}
\newcommand{\cov}{\mathrm{cov}}
\newcommand{\eff}{\mathrm{eff}}
\newcommand{\Ho}{\mathrm{Ho}}
\newcommand{\id}{\mathrm{id}}

\newcommand{\KGL}{\mathrm{KGL}}
\newcommand{\MM}{\mathrm{M}}
\newcommand{\MGL}{\mathrm{MGL}}
\newcommand{\modd}{\mathrm{mod}}
\newcommand{\MU}{\mathrm{MU}}

\newcommand{\pre}{\mathrm{pre}}
\newcommand{\Pre}{\mathrm{Pre}}

\newcommand{\pur}{\mathfrak{pur}} %

\newcommand{\SH}{\mathrm{SH}}
\newcommand{\Sm}{\mathrm{Sm}}
\newcommand{\Sp}{\mathrm{Sp}}
\newcommand{\Spc}{\mathrm{Spaces}}
\newcommand{\Spt}{\mathrm{Spectra}}
\newcommand{\Thom}{\mathrm{Th}} %
\newcommand{\tr}{\mathrm{tr}}
\newcommand{\unst}{\mathrm{H}} %

\newcommand{\exch}[1]{Ex(#1_*^*,\sm)} %
\newcommand{\exmodl}{Ex_{\otimes}^{*!}} %
\newcommand{\exmodr}{Ex_{\otimes}^{!*}} %
\newcommand{\seq}{\alpha} %
\newcommand{\steen}{\mathcal{A}} %

\DeclareMathOperator{\cha}{char} %

\DeclareMathOperator{\Frac}{Frac}
\DeclareMathOperator*{\hocolim}{hocolim}
\DeclareMathOperator*{\holim}{holim}
\DeclareMathOperator{\Hom}{Hom}

\DeclareMathOperator{\Pic}{Pic}
\DeclareMathOperator{\Spec}{Spec}
\DeclareMathOperator{\Tor}{Tor}

\newcommand{\bigcol}{3.6pc}

\begin{document}

\title{Towards the dual motivic Steenrod algebra in positive characteristic} 
\date{\today}

\author{Martin Frankland}     
\address{University of Regina}        
\email{Martin.Frankland@uregina.ca}

\author{Markus Spitzweck}
\address{Universit\"at Osnabr\"uck}       
\email{markus.spitzweck@uni-osnabrueck.de}

\begin{abstract}
The dual motivic Steenrod algebra with mod $\ell$ coefficients was computed by Voevodsky over a base field of characteristic zero, and by Hoyois, Kelly, and {\O}stv{\ae}r over a base field of characteristic $p \neq \ell$. In the case $p = \ell$, we show that the conjectured answer is a retract of the actual answer. We also describe the slices of the algebraic cobordism spectrum $\MGL$: we show that the conjectured form of $s_n \MGL$ is a retract of the actual answer.
\end{abstract}



\maketitle

\tableofcontents

\section{Introduction}

In \cite{Voevodsky02ope}, Voevodsky proposed a list of major conjectures in motivic homotopy theory, some of which were known over a base field of characteristic zero. In this paper, we prove a weak form of two of the Voevodsky Conjectures over a base field of positive characteristic $p$: one about the mod $p$ dual motivic Steenrod algebra, and one about the slices of the algebraic cobordism spectrum $\MGL$.
The computation of the mod $p$ topological Steenrod algebra (or its dual) was a milestone in algebraic topology,
both as a computational and a conceptual tool.
The motivic Steenrod algebra over a field of characteristic $0$ played a major role in Voevodsky's celebrated proofs of the Milnor and motivic Bloch--Kato Conjectures,
e.g., there had to be enough operations to define motivic Margolis cohomology. Also for the Bloch--Kato Conjecture the exact structure of the motivic Steenrod algebra
had to be determined \cite{Voevodsky10}. %
This latter point is missing so far for the mod $p$ motivic Steenrod algebra over fields of characteristic $p$. 
Mod $p$ motivic cohomology over fields of characteristic $p$, the coefficients of the sought after dual motivic Steenrod algebra, is understood.
Only the Milnor-line inside the bigraded motivic cohomology sheaves 
is non-zero, a result of Geisser and Levine \cite{GeisserL00}, and it is given by the logarithmic de Rham sheaves (the Bloch--Gabber--Kato theorem \cite{BlochK86}).

\subsection*{Dual motivic Steenrod algebra}

Let $\MM\F_{\ell}$ denote the mod $\ell$ motivic Eilenberg--MacLane spectrum. There are certain classes in the dual motivic Steenrod algebra $\pi_{*,*} (\MM\F_{\ell} \sm \MM\F_{\ell})$, the ``Milnor basis'', which together induce a map of $\MM\F_{\ell}$-modules 
\begin{equation}\label{eq:DecompoIntro}
\xymatrix{
\bigoplus_{\seq} \Si^{p_{\seq},q_{\seq}} \MM\F_{\ell} \ar[r]^-{\Psi} & \MM\F_{\ell} \sm \MM\F_{\ell}. \\
}
\end{equation}
The map $\Psi$ is conjectured to be an equivalence \cite{Voevodsky02ope}*{Conjecture~4} over any base scheme. 
Voevodsky proved this conjecture when the base scheme $S$ is a field of characteristic zero \cite{Voevodsky03red}. Hoyois, Kelly, and {\O}stv{\ae}r proved it when $S$ is essentially smooth over a field of characteristic $p \neq \ell$ \cite{HoyoisKO17}.

Our first main result is the following (see Theorem~\ref{thm:RetractSteenrod}).

\begin{theorem}\label{thm:RetractSteenrodIntro}
Assume that the base scheme $S$ is essentially smooth over a field of characteristic $p > 0$. Then the map of $\MM\F_p$-modules %
$\Psi$ in \eqref{eq:DecompoIntro} is a split monomorphism, i.e., admits a retraction. 
\end{theorem}

It was kindly pointed out to us by Marc Hoyois that Theorem~\ref{thm:RetractSteenrodIntro} could be used to produce new motivic cohomology operations in positive characteristic. 
In the meantime, this was worked out by Primozic \cite{Primozic20}.

\subsection*{Slices of algebraic cobordism}

The algebraic cobordism spectrum $\MGL$ comes with a canonical ring map $\la \colon L_{*} \to \pi_{2*,*} \MGL$ from the Lazard ring $L \cong \Z[x_1, x_2, \ldots]$. Here we use the grading $\abs{x_i}=i$, instead of the grading $\abs{x_i} = 2i$ which is customary in topology \cite{Ravenel04}*{\S A2.1}. %
Over a field, the map $\la$ induces in turn a canonical map of $\MM\Z$-modules %
\begin{equation}\label{eq:ApproxSliceDMIntro}
\xymatrix{
\Si^{2n,n} \MM\Z \ot_{\Z} L_{n} \ar[r]^-{\psi} & s_n \MGL \\
}
\end{equation}
to the $n$\textsuperscript{th} slice $s_n \MGL$. Voevodsky conjectured that $\psi$ is an equivalence over any base field %
\cite{Voevodsky02ope}*{Conjecture~5}. Over a field of characteristic zero, the conjecture was proved in \cite{Spitzweck10} as an application of the Hopkins--Morel isomorphism. 
Over a base field of characteristic $p$, Hoyois showed that $\psi$ induces an equivalence after inverting $p$ \cite{Hoyois15}*{\S 8.3}.

Here is our second main result (see Theorem~\ref{thm:RetractSlice}).

\begin{theorem}\label{thm:RetractSliceIntro}
Assume that the base scheme $S$ is essentially smooth over a field of characteristic $p > 0$. Then the map of $\MM\Z$-modules $\psi$ in \eqref{eq:ApproxSliceDMIntro} is a split monomorphism, i.e., admits a retraction.
\end{theorem}

\subsection*{A possible strategy}

As a reduction step towards proving Voevodsky's Conjecture~4, we show the following conditional statement (see Proposition~\ref{pr:ReduceLisse}). Let $\hatt{M}\Z$ denote the version of the Eilenberg--MacLane spectrum constructed in \cite{Spitzweck18}.

\begin{proposition}\label{pr:ReduceLisseIntro}
Let $\Z_p$ denote the $p$-adic integers. \emph{If} the object $\hatt{M}\F_{p}$ of $\SH(\Z_p)$ lies in the localizing triangulated subcategory of $\SH(\Z_p)$ generated by the strongly dualizable objects, then the map of $\MM\F_p$-modules $\Psi$ in Theorem~\ref{thm:RetractSteenrodIntro} is an equivalence.
\end{proposition}

\subsection*{Relation to the Hopkins--Morel--Hoyois isomorphism}

For $n \geq 1$, let $a_n \in \pi_{2n,n} \MGL$ denote classes coming from appropriate generators of the Lazard ring $L$; see \cite{Hoyois15} for details. 
This induces a canonical map in $\SH(S)$
\[
\Phi \colon \MGL / (a_1, a_2, \ldots) \to \MM\Z.
\]

Hoyois proved the following \cite{Hoyois15}; the case $\cha \kk = 0$ is due to Hopkins and Morel.

\begin{theorem}
Assume that the base scheme $S$ is essentially smooth over a field $\kk$ of characteristic exponent $p$ (i.e., set $p=1$ in the case $\cha \kk = 0$). 
Then $\Phi$ induces an isomorphism
\[
\xymatrix @C=\bigcol {
\MGL / (a_1, a_2, \ldots) \left[ 1/p \right] \ar[r]^-{\Phi [1/p]}_-{\cong} & \MM\Z \left[ 1/p \right]. \\
}
\]
\end{theorem}
The theorem was further generalized to $S$ being essentially smooth over a regular $\kk$-scheme \cite{CisinskiD15}*{Remark~3.7}.

\begin{conjecture}\label{conj:HopkinsMorel}
The map $\Phi$ %
is an isomorphism for any base scheme $S$.
\end{conjecture}

The work in \cite{Hoyois15} reduces the problem as follows.

\begin{theorem}
Over a field of characteristic $p > 0$, if $\MM\F_p \sm \Phi$ is an isomorphism, then so is $\Phi$ itself.
\end{theorem}

The calculation of the dual Steenrod algebra in \cite{HoyoisKO17} was the main ingredient in the proof of the Hopkins--Morel--Hoyois isomorphism \cite{Hoyois15}. 
If one could prove the conjecture on the dual Steenrod algebra in the case $\ell = p$, one should be able to adapt the argument in \cite{Hoyois15} to prove that $\Phi$ is an isomorphism, without inverting $p$.

\subsection*{Organization}

In Section~\ref{sec:Setup}, we describe the setup of motivic spectra and the six functor formalism.  In Section~\ref{sec:DualSteenrod}, we recall facts about the dual motivic Steenrod algebra and Voevodsky's Conjecture~4. Section~\ref{sec:DVR} introduces a discrete valuation ring $D$ of mixed characteristic, which we use as a bridge between characteristic $p$ and characteristic $0$. Our main technical ingredient is Lemma~\ref{lem:EMopenClosed}. Our main results (Theorem~\ref{thm:RetractSteenrodIntro}, Theorem~\ref{thm:RetractSliceIntro}, and Proposition~\ref{pr:ReduceLisseIntro}) are then proved in Sections~\ref{sec:RetractSteenrod}, \ref{sec:RetractSlices}, and \ref{sec:ReduceLisse} respectively. 
Some technical points are deferred to %
Section~\ref{sec:Purity}, notably a purity result for motivic cohomology (Proposition~\ref{pr:ShriekHZ}).

\subsection*{Acknowledgments}
We thank Tobias Barthel, David Carchedi, Elden Elmanto, Marc Hoyois, and Beren Sanders for useful discussions. We also thank Dan Isaksen for helpful comments on an earlier version, as well as Fr\'ed\'eric D\'eglise and Adeel Khan for sharing their preliminary work. 
Moreover, we thank the anonymous referee for their careful reading and their suggestions.

This project was supported by a grant of the Deutsche Forschungsgemeinschaft SPP 1786: Homotopy Theory and Algebraic Geometry. Frankland acknowledges the support of the Natural Sciences and Engineering Research Council of Canada (NSERC), Discovery Grant RGPIN-2019-06082. 

\section{Setup, notations, and conventions}\label{sec:Setup}

\subsection{Motivic spaces and spectra}

The $\A^1$-homotopy theory of schemes was developed by Morel and Voevodsky \cite{MorelV99}. See \cite{HoyoisKO17}*{\S 2.1}, \cite{Hoyois15}*{\S 2}, or \cite{Voevodsky10}*{\S 1.1} for background material.

\begin{notation}
We will work over schemes that are Noetherian, separated, and of finite Krull dimension. 
Such a scheme $S$ is called a \Def{base scheme}. 
We will often assume that $S$ is over some field $\kk$ of characteristic $p$. 
Let $\base$ denote the category of base schemes, viewed as a full subcategory of all schemes. 
\end{notation}

Let $\Sm_S$ denote the category of separated smooth schemes of finite type over $S$. Under the assumptions on $S$, the category $\Sm_S$ is essentially small. The category of \Def{motivic spaces} over $S$ is
\[
\Spc(S) \dfn L_{\A^1} L_{\mathrm{Nis}} s\Pre(\Sm_S),
\]
where the category $s\Pre(\Sm_S)$ of simplicial presheaves is endowed with the projective model structure. Here, $L_{\mathrm{Nis}}$ denotes the (left) Bousfield localization with respect to Nisnevich hypercovers, and $L_{\A^1}$ is the Bousfield localization with respect to the projection maps $X \x \A^1 \to X$, where $X$ ranges over representatives of isomorphism classes of objects in $\Sm_S$.

The category of \Def{motivic spectra} over $S$ is the category of symmetric $\PP^1$-spectra in pointed motivic spaces:
\[
\Spt(S) \dfn \Sp^{\Si}(\Spc_*(S), \PP^1). 
\]
Various model structures on motivic spectra are discussed in \cite{Hovey01} and \cite{Jardine00}.

\begin{notation}\label{nota:Fibered}
Over a base scheme $S$, denote by:
\begin{itemize}
\item $\unst(S)$ the motivic unstable homotopy category $\Ho (\Spc(S))$.
\item $\SH(S)$ the motivic stable homotopy category $\Ho (\Spt(S))$.
\end{itemize}
Under the assumptions on $S$, $\SH(S)$ is a compactly generated triangulated category with arbitrary coproducts.
\end{notation}

\begin{notation}
In $\SH(S)$, denote the bigraded motivic sphere spectra by
\begin{align*}
S^{p,q} &= (S^1)^{\sm (p-q)} \sm_S \G_m^{\sm q} \\
&\dfn (\Si^{\infty} S^1)^{\sm (p-q)} \sm_S (\Si^{\infty} \G_m)^{\sm q}
\end{align*}
and the corresponding suspension functors by
\[
\Si^{p,q} X = S^{p,q} \sm_S X = X(q)[p]. 
\]
Denote the sphere spectrum by $\unit_S = S^{0,0} = \Si^{\infty} S_+$, which is the unit for the smash product $\sm_S$. 
\end{notation}

Recall the following convention from \cite{Spitzweck18}*{\S 2}. An $E_{\infty}$ structure is understood with respect to (the image of) the linear isometries operad (see \cite{EKMM97}*{\S I.3}). To be more precise, taking normalized chains of the linear isometries operad (which is an operad in spaces) yields an operad $E$ in $\Ch_{\geq 0}(\Ab)$, the category of non-negatively graded chain complexes of abelian groups. Our algebras will live in categories $\cat{C}$ which receive a symmetric monoidal 
functor $s$ from $\Ch_{\geq 0}(\Ab)$, and then we will call an algebra over the operad $s(E)$ an $E_{\infty}$ algebra in $\cat{C}$. 
In particular, an $E_{\infty}$ motivic ring spectrum (over a base scheme $S$) means an $E_{\infty}$ algebra in the category of motivic spectra $\Spt(S)$.

\begin{notation}
For an $E_{\infty}$ motivic ring spectrum $R_S$ over $S$, denote by $\hMod{R_S}$ the homotopy category of (left) $R_S$-modules. It comes with a forgetful functor $U_S \colon \hMod{R_S} \to \SH(S)$, as well as a (derived) smash product $\sm_{R_S}$ 
that makes $\hMod{R_S}$ symmetric monoidal. 
In this notation, $\SH(S) = \hMod{\unit_S}$. 
\end{notation}

\subsection{Dependence on the base scheme}

By varying the base scheme $S$, $\SH$ defines a category fibered over $\base$. %
Ayoub showed that $\SH$ satisfies the \emph{six functor formalism} of Grothendieck \cite{Ayoub07I}*{\S 1.4.1}. See also %
\cite{CisinskiD19}*{\S 1} 
and \cite{Hoyois14}*{\S 2} for more details.

\begin{notation}
Let $f \colon S \to T$ be a map of schemes. Denote by $f_* := Rf_*$ the total right derived direct image functor
\[
f_* \colon \SH(S) \to \SH(T)
\]
induced on the homotopy categories of motivic spectra. Denote by $f^* := Lf^*$ the total left derived inverse image functor 
\[
f^* \colon \SH(T) \to \SH(S)
\]
which is left adjoint to $f_*$.

If moreover the map $f$ is separated and of finite type, then it induces a further (derived) adjunction $f_! \colon \SH(S) \rla \SH(T) \colon f^!$. There is a natural transformation $f_! \to f_*$, which is an isomorphism if $f$ is proper. In particular, a closed immersion $i$ satisfies $i_! = i_*$. An open immersion $j$ satisfies $j^! = j^*$. %
\end{notation}

The functor $f^*$ is always strong monoidal, i.e., there are natural isomorphisms \linebreak[4]
$f^*A \sm_S f^*B \ral{\cong} f^*(A \sm_T B)$ and $\unit_{S} \ral{\cong} f^*(\unit_T)$, 
cf.\ %
\cite{CisinskiD19}*{Definitions~1.1.21, 1.4.2, Example~1.4.3} 
and \cite{Hoyois14}*{\S 2}. 
The right adjoint $f_*$ is lax monoidal, though need not be %
strong monoidal. For a closed immersion $i \colon Z \to X$, $i_*$ does preserve the smash product, but need not preserve the tensor unit $\unit$ \cite{Ayoub07I}*{Lemme~2.3.6, Remarque~2.3.7}. More precisely, the comparison map $\unit_X \to i_*(\unit_Z)$ is an isomorphism if and only if the open complement $X \setminus Z$ is empty \cite{CisinskiD19}*{\S 2.3.d}.

The map of base schemes $f \colon S \to T$ induces an inverse image functor $f^*$ on module categories, compatible with the functor $f^*$ on underlying spectra. More precisely, let $R_T$ be an $E_{\infty}$ ring spectrum in $\Spt(T)$. Then the following diagram commutes up to natural isomorphism:
\begin{equation}\label{eq:PullbackForget}
\xymatrix{
\hMod{R_T} \ar[d]_{U_T} \ar[r]^-{f^*} & \hMod{f^* R_T} \ar[d]^{U_S} \\
\SH(T) \ar[r]^-{f^*} & \SH(S), \\
}
\end{equation}
where the downward arrows denote the forgetful functors. 
See also \cite{CisinskiD15}*{Proposition~3.10, Lemma~3.20}.

The motivic Eilenberg--MacLane spectrum $\MM\Z = \{ \MM\Z_S \}_S$ is defined over an arbitrary base scheme $S$. Let us recall the structure relating $\MM\Z_S$ and $\MM\Z_T$ over different base schemes. %

\begin{definition}
A \Def{$\base$-fibered} motivic spectrum $X$ is a section of the projection functor $\SH \to \base$ defining $\SH$ as a fibered category over all base schemes. Explicitly, $X = \{ X_S \}$ consists of:
\begin{itemize}
\item For every base scheme $S$, an object $X_S$ of $\SH(S)$.
\item For every map of base schemes $f \colon S \to T$, a %
map $f^* X_T \to X_S$ in $\SH(S)$.
\item For every pair of composable maps $S \ral{f} T \ral{g} U$ of base schemes, the diagram in $\SH(S)$
\begin{equation}\label{eq:CompoCompar}
\xymatrix{
f^* g^* X_U \ar[d]_{\cong} \ar[r] & f^* X_T \ar[r] & X_S \\
(gf)^* X_U \ar[urr] & & \\
}
\end{equation}
commutes, where $f^* g^* \ral{\cong} (gf)^*$ 
is the connection isomorphism for the $2$-functor $\SH \colon \base \to \Cat$.
\end{itemize}
\end{definition}

\begin{definition}\label{def:DM}
A \Def{$\base$-fibered} motivic $E_{\infty}$ ring spectrum $R$ is a section of the fibered category $\Ho (E_{\infty}-\Spt)$.

As explained in %
\cite{CisinskiD19}*{\S 7.2}, 
for a map of base schemes $f \colon S \to T$ and an $R_T$-module $M_T$, the pullback $f^* M_T$ naturally becomes an $f^* R_T$-module, which in turn can be extended to an $R_S$-module via the structure map $f^* R_T \to R_S$, that is:
\[
f^*_{\modd}(M_T) \dfn R_S \sm_{f^* R_T} f^* M_T.
\]
This construction defines the pullback functors on modules $f^*_{\modd} \colon \hMod{R_T} \to \hMod{R_S}$ making $\hMod{R}$ into a $\base$-fibered category.

We say that $R$ is \Def{Cartesian} if the structure maps $f^* R_T \ral{\cong} R_S$ are isomorphisms. 
(This property of $R$ is also called \emph{absolute} in the literature \cite{Deglise18ori}*{Definition~1.2.1}.) 
In that case, the pullback in modules $f^*_{\modd}(M_T)$ has the same underlying spectrum as the pullback in spectra $f^*(M_T)$.
\end{definition}

The functor $f^!$ has structure maps $\exmodr \colon f^!A \sm f^*B \to f^!(A \sm B)$ and $\exmodl \colon f^*B \sm f^!A \to f^!(B \sm A)$ satisfying certain equations, making $f^!$ an \emph{$f^*$-bimodule} in the terminology of \cite{Ayoub07I}*{\S 2.3.5}; cf.\ Notation~\ref{nota:ExchangeModule}. 
Given an $R_T$-module $M_T$, $f^!(M_T)$ canonically becomes an $f^*(R_T)$-module via the structure map $\exmodl$. If the ring spectrum $R$ is Cartesian, then this construction defines the functor on modules \linebreak $f^!_{\modd} \colon \hMod{R_T} \to \hMod{R_S}$ which is part of the six functor formalism for $\hMod{R}$. In this case, the underlying spectrum of $f^!_{\modd}(M_T)$ is $f^!(M_T)$.

\subsection{Choice of models}

Our main tool in this paper is the version of the Eilenberg--MacLane spectrum $\hatt{M}\Z$ constructed in \cite{Spitzweck18}. It is constructed as an $E_{\infty}$ algebra in symmetric $\Z[\G_m, \{1\}]$-spectra %
of complexes of sheaves of abelian groups on $\Sm_S$. To obtain well-behaved symmetric monoidal homotopy categories of $\hatt{M}\Z$-modules as described above, one can use the strictification results of \cite{Hornbostel13}*{\S 3} to replace $\hatt{M}\Z$ by a (strictly) commutative monoid in some symmetric monoidal model category of motivic spectra. %
Alternately, one can use the results of \cite{Robalo15}*{\S 2.4} to construct the desired homotopy categories by working only at the $\infty$-categorical level.

In this paper, our statements will take place in homotopy categories, such as $\SH(S)$ or $\hMod{R_S}$. To emphasize this, we will speak of isomorphisms instead of equivalences. 

We will consider some homotopy colimits of sequences $X_0 \to X_1 \to \cdots$ in those homotopy categories. Though our examples lift to sequences in model categories, recall that a sequential homotopy colimit can be constructed (non functorially) within the triangulated category \cite{Neeman01}*{Definition~1.6.4}.

Also, to ease the notation, we may omit the last map of an exact triangle, i.e., write $A \to B \to C$ instead of $A \to B \to C \to \Si^{1,0} A$.

\section{The dual Steenrod algebra}\label{sec:DualSteenrod}

Background on motivic Eilenberg--MacLane spectra can be found in \cite{Voevodsky10}*{\S 3.2}, \cite{RondigsO08}*{\S 2.4}, \cite{HoyoisKO17}*{\S 2.3}, or \cite{Hoyois15}*{\S 4.2}. 

\begin{notation}\label{nota:EMspectra}
For a base scheme $S$, %
let $\MM \Z_S$ denote the motivic Eilenberg--MacLane spectrum in $\SH(S)$ constructed by Voevodsky. %
There is a version with coefficients in an abelian group $A$, which we denote
$\MM A_S = \MM\Z_S \ot_{\Z} A$.
If $S$ is an affine scheme $\Spec(R)$, we drop $\Spec$ from the notation and write $\MM A_R$. %

We write $\MM\F_{\ell}^S$ instead of $\MM\F_{\ell,S}$ or $M(\F_{\ell})_S$ to avoid overcrowding the subscript.
We drop the base scheme $S$ from the notation when there is no risk of confusion.
\end{notation}

We recall the structure of the dual motivic Steenrod algebra, described for instance in \cite{Voevodsky03red}*{\S 12} \cite{Hoyois15}*{\S 5.3}, 
\cite{Spitzweck18}*{Theorem~10.26}. 
For the classical analogue in topology, see \cite{Milnor58}*{\S 5}. There are certain classes $\tau_i, \xi_i \in \pi_{*,*} (\MM\F_{\ell} \sm \MM\F_{\ell})$ with bidegrees
\begin{align*}
\abs{\tau_i} &= (2 {\ell}^i-1, {\ell}^i-1), i \geq 0 \\
\abs{\xi_i} &= (2 {\ell}^i-2, {\ell}^i-1), i \geq 1.
\end{align*}
Consider sequences $\seq = (\ep_0, r_1, \ep_1, r_2, \ep_2, \ldots)$ with $\ep_i \in \{0,1\}$, $r_i \geq 0$, and only finitely many non-zero terms. To each such sequence, one assigns the monomial
\[
\om(\seq) \dfn \tau_0^{\ep_0} \xi_1^{r_1} \tau_1^{\ep_1} \xi_2^{r_2} \tau_2^{\ep_2} \cdots \in \pi_{*,*} (\MM\F_{\ell} \sm \MM\F_{\ell})
\]
of bidegree $\abs{\om(\seq)} = (p_{\seq}, q_{\seq})$. Together, these monomials form the motivic analogue of the Milnor basis. 
Consider the induced map of $\MM\F_{\ell}$-modules 
\begin{equation}\label{eq:Decompo}
\xymatrix{
\bigoplus\limits_{\seq} \Si^{p_{\seq},q_{\seq}} \MM\F_{\ell} \ar[r]^-{\Psi^S} & \MM\F_{\ell} \sm \MM\F_{\ell}. \\
}
\end{equation}
Here we view $\MM\F_{\ell} \sm \MM\F_{\ell}$ as a left $\MM\F_{\ell}$-module by acting on the left smash factor. 

\begin{theorem}\label{thm:DualSteenrod}
Assume that the base scheme $S$ is essentially smooth over a field of characteristic $p$. 
Then for $\ell \neq p$, the map %
$\Psi^S$ in $\hMod{\MM\F_{\ell}^S}$ displayed in Equation~\eqref{eq:Decompo} is an isomorphism.%
\end{theorem}

For $S$ a field of characteristic zero, this result is due to Voevodsky 
(using \cite{Voevodsky03red}*{Corollary~12.5}; cf.\ \cite{DuggerI10}*{Lemma~6.3}, \cite{Hoyois15}*{\S 5.4}). 
It was generalized to the stated form in \cite{HoyoisKO17}*{Theorem~1.1, Corollary~3.4}. It is conjectured to hold for more general base schemes \cite{Voevodsky02ope}*{Conjecture~4}.

\subsection{Eilenberg--MacLane spectra and base change}

Recall Voevodsky's Conjecture~17 from \cite{Voevodsky02ope}.

\begin{conjecture}\label{conj:PullbackHZ}
Let $f \colon S \to T$ be a map of base schemes. Then the induced map in $\SH(S)$
\[
f^* \MM\Z_T \to \MM\Z_S
\]
is an isomorphism.
\end{conjecture}

For the record, we recall the following:

\begin{definition}
\cite{Hoyois15}*{Definition~A.1} A map of base schemes $f \colon S \to T$ is \Def{essentially smooth} if $S$ is a cofiltered limit $S = \lim_{\al} S_{\al}$ of smooth $T$-schemes where the transition maps $S_{\be} \to S_{\al}$ in the diagram are affine and dominant.
\end{definition}

\begin{theorem}\label{thm:PullbackHZ}
Let $f \colon S \to T$ be a map of base schemes.
\begin{enumerate}
\item \label{item:PullbackSmooth} If $f$ is essentially smooth, then the map $f^* \MM\Z_{T} \ral{\cong} \MM\Z_S$ is an isomorphism \cite{HoyoisKO17}*{Theorem~2.12} \cite{Hoyois15}*{Theorem~4.18}.
\item If $S$ and $T$ are regular over a base field of characteristic $p$, then the map \linebreak $f^* \MM\Z[1/p]_T \ral{\cong} \MM\Z[1/p]_S$ is an isomorphism \cite{CisinskiD15}*{Corollary~3.6}. 
\end{enumerate}
\end{theorem}

\subsection{The Milnor basis and base change}

\begin{definition}
Let $R$ be a motivic $E_{\infty}$ ring spectrum.

\begin{enumerate}
\item A left \Def{$R$-module} $M$ is a section of the fibered category $\hMod{R}$. Explicitly, $M = \{ M_S \}$ consists of:
 \begin{itemize}
 \item For every base scheme $S$, a left $R_S$-module $M_S$ in the homotopy category $\hMod{R_S}$. We denote the action map %
 by $\mu_S \colon R_S \sm M_S \to M_S$.
 \item For every map of base schemes $f \colon S \to T$, %
a map $f^* M_T \to M_S$ in $\hMod{f^* R_T}$; equivalently, a map in $\SH(S)$ making the following diagram in $\SH(S)$ commute:
\[
\xymatrix @R-0.5pc {
f^* (R_T \sm M_T) \ar[r]^-{f^* \mu_T} & f^* M_T \ar[dd] \\
f^* R_T \sm f^* M_T \ar[d] \ar[u]^{\cong} & \\
R_S \sm M_S \ar[r]^{\mu_S} & M_S. \\
}
\]
As usual, the coherence condition~\eqref{eq:CompoCompar} %
holds.
 \end{itemize}
\item A class $\xi \in \pi_{p,q} M$ means a family $\xi = \{ \xi_S \}$ with $\xi_S \in \pi_{p,q} M_S$ for every base scheme $S$. Such a class $\xi$ is \Def{compatible with base change} if for every map of base schemes $f \colon S \to T$, the induced map
\[
\xymatrix{
\pi_{p,q} M_T \ar[r]^-{f^*} & \pi_{p,q} (f^* M_T) \ar[r]^{} & \pi_{p,q} M_S
}
\]
sends $\xi_T$ to $\xi_S$. 
In other words, the following diagram in $\SH(S)$ commutes:
\[
\xymatrix{
f^* (\Si^{p,q} \unit_T) \ar[r]^-{f^* \xi_T} & f^* M_T \ar[d] \\
\Si^{p,q} \unit_S \ar[r]^-{\xi_S} \ar[u]^{\cong} & M_S. \\
}
\]
\end{enumerate}
\end{definition}

\begin{lemma}\label{lem:CartesianClasses}
Let $R$ be a motivic $E_{\infty}$ ring spectrum, $M$ an $R$-module, and $\xi \in \pi_{p,q} M$ a class which is compatible with base change. Then its adjunct map of $R$-modules $\tild{\xi} \colon \Si^{p,q} R \to M$ %
is also compatible with base change, in the sense that for every map of base schemes $f \colon S \to T$, the following diagram in $\hMod{f^* R_T}$ commutes: %
\[
\xymatrix @R-0.5pc {
f^* (\Si^{p,q} R_T) \ar[r]^-{f^* \tild{\xi}_T} & f^* M_T \ar[dd] \\
\Si^{p,q} f^* R_T \ar[d] \ar[u]^{\cong} & \\
\Si^{p,q} R_S \ar[r]^-{\tild{\xi}_S} & M_S. \\
}
\]
\end{lemma}

\begin{proof}
This follows from a straightforward diagram chase.
\end{proof}

\begin{lemma}\label{lem:MilnorCartesian}
The classes $\tau_i, \xi_i \in \pi_{*,*} (\MM\F_{\ell} \sm \MM\F_{\ell})$ are compatible with base change.
\end{lemma}

\begin{proof}
Write $\steen^S_{*,*} \dfn \pi_{*,*} (\MM\F_{\ell}^S \sm \MM\F_{\ell}^S)$ for short. Let us recall the construction of the classes $\tau_i, \xi_i \in \steen_{*,*}$ from %
\cite{Spitzweck18}*{\S 10.2}.

Let $W_{S,n}$ denote the $\G_m$-torsor over $\PP^{\infty}_S$ corresponding to the line bundle $\mathcal{O}_{\PP^{\infty}_S}(-n)$, with projection map $W_{S,n} \to \PP^{\infty}_S$, viewed as a map of motivic spaces over $S$. 
The motivic cohomology %
of projective space $\PP^{\infty}_S$ is
\[
H^{*,*}(\PP^{\infty}_S) = H^{*,*}(S) \Brack{\si}
\]
where $\si \in H^{2,1}(\PP^{\infty}_S)$ is the class of the line bundle $\mathcal{O}_{\PP^{\infty}_S}(-1)$. 
The cohomology groups $H^{*,*}(W_{S,n})$ sit in a long exact sequence 
\[
\xymatrix @C-0.6pc {
\cdots \ar[r] & H^{*-2,*-1}(S) \Brack{\si} \ar[r]^-{n \si} & H^{*,*}(S) \Brack{\si} \ar[r] & H^{*,*}(W_{S,n}) \ar[r] & H^{*-1,-1}(S) \Brack{\si} \ar[r] & \cdots \\
}
\]
from which $H^{*,*}(W_{S,n})$ can be computed \cite{Spitzweck18}*{Theorem~10.16}. Let $c_S \in H^{2,1}(\PP^{\infty}_S; \F_{\ell})$ denote the image of $\si \in H^{2,1}(\PP^{\infty}_S)$ under reduction of coefficients modulo $\ell$, and 
let $v_S \in H^{2,1}(W_{S,\ell}; \F_{\ell})$ denote the restriction of $c_S$ along the map $W_{S,\ell} \to \PP^{\infty}_S$. 
The canonical mod-$\ell$ orientation of $\MM\F_{\ell}$ consists of the classes
\[
c_S \in H^{2,1}(\PP^{\infty}_S; \F_{\ell}) \quad \text{and} \quad u_S \in H^{1,1}(W_{S,\ell}; \F_{\ell})
\]
where $u_S$ is related to $v_S$ via an exact sequence calculation \cite{Spitzweck18}*{Lemma~10.6}. 
The classes $\tau_i$ and $\xi_i$ are defined as the coefficients appearing in the coaction:
\[
\xymatrix @R-1.2pc {
H^{*,*}(W_{S,\ell}; \F_{\ell}) \ar[r]^-{\coact} & \steen^S_{-*,-*} \hat{\ot}_{\pi_{-*,-*} \MM\F_{\ell}^S} H^{*,*} (W_{S,\ell}; \F_{\ell}) \\
u \ar@{|->}[r] & u + \sum_{i \geq 0} \tau_i \ot v^{\ell^i} \\
v \ar@{|->}[r] & v + \sum_{i \geq 1} \xi_i \ot v^{\ell^i}. \\
}
\]
Now the claim follows from these two facts. First, the coaction is compatible with base change, i.e., the following diagram of abelian groups commutes:
\[
\xymatrix @R-0.4pc {
H^{p,q}(W_{T,\ell}; \F_{\ell}) \ar[d] \ar[r]^-{\coact_T} & \steen^T_{-*,-*} \hat{\ot}_{\pi_{-*,-*} \MM\F_{\ell}^T} H^{*,*} (W_{T,\ell}; \F_{\ell}) \ar[d] \\
H^{p,q}(f^* W_{T,\ell}; \F_{\ell}) & \steen^S_{-*,-*} \hat{\ot}_{\pi_{-*,-*} \MM\F_{\ell}^S} H^{*,*} (f^* W_{T,\ell}; \F_{\ell}) \\
H^{p,q}(W_{S,\ell}; \F_{\ell}) \ar[r]^-{\coact_S} \ar[u]^{\cong} & \steen^S_{-*,-*} \hat{\ot}_{\pi_{-*,-*} \MM\F_{\ell}^S} H^{*,*} (W_{S,\ell}; \F_{\ell}). \ar[u]_{\cong} \\
}
\]
The downward maps in the top part are induced by the functor $f^*$ and by the map $f^* \MM\F_{\ell}^T \to \MM\F_{\ell}^S$. The upward maps in the bottom part are induced by the map $f^* W_{T,\ell} \ral{\simeq} W_{S,\ell}$, which is an isomorphism in $\unst(S)$.

Second, the classes $u_{S} \in H^{1,1}(W_{S,\ell}; \F_{\ell})$ and $v_{S} \in H^{2,1}(W_{S,\ell}; \F_{\ell})$ are compatible with base change, i.e., the composite
\[
\xymatrix @R-1.2pc {
H^{*,*}(W_{T,\ell}; \F_{\ell}) \ar[r] & H^{*,*}(f^* W_{T,\ell}; \F_{\ell}) & H^{*,*}(W_{S,\ell}; \F_{\ell}) \ar[l]_-{\cong} \\
u_{T} \ar@{|->}[rr] & & u_{S} \\
v_{T} \ar@{|->}[rr] & & v_{S} \\ 
}
\]
sends the classes $u$ and $v$ as indicated.
\end{proof}

The decomposition map $\Psi^S$ is compatible with base change in the following sense.

\begin{lemma}\label{lem:CompatDecompo}
Let $f \colon S \to T$ be a map of base schemes. Then the following diagram in $\hMod{\MM\F_{\ell}^S}$ commutes:
\[
\xymatrix @R-0.5pc {
f^* \left( \bigop_{\seq} \Si^{p_{\seq},q_{\seq}} \MM\F_{\ell}^{T} \right) \ar[r]^-{f^* \Psi^T} & f^* \left( \MM\F_{\ell}^{T} \sm_{T} \MM\F_{\ell}^{T} \right) \\
\bigop_{\seq} \Si^{p_{\seq},q_{\seq}} f^* \MM\F_{\ell}^{T} \ar[d] \ar[u]^{\cong} &  f^* \MM\F_{\ell}^{T} \sm_{S} f^* \MM\F_{\ell}^{T} \ar[u]_{\cong} \ar[d] \\
\bigop_{\seq} \Si^{p_{\seq},q_{\seq}} \MM\F_{\ell}^{S} \ar[r]^-{\Psi^S} &  \MM\F_{\ell}^{S} \sm_{S} \MM\F_{\ell}^{S}. \\
}
\]
\end{lemma}

\begin{proof}
The algebra structure of the dual Steenrod algebra is compatible with pullbacks, i.e., the induced map
\[
\xymatrix{
\pi_{*,*} (\MM\F_{\ell}^T \sm \MM\F_{\ell}^T) \ar[r] & \pi_{*,*} (\MM\F_{\ell}^S \sm \MM\F_{\ell}^S)
}
\]
is a map of algebras. This follows from the fact that $f^*$ is strong monoidal and the map $f^* \MM\F_{\ell}^T \to \MM\F_{\ell}^S$ is a map of $E_{\infty}$ ring spectra. %
Now the statement follows from Lemmas~\ref{lem:MilnorCartesian} and \ref{lem:CartesianClasses}.
\end{proof}

\begin{remark}
For an essentially smooth map $f$ over a field, a similar statement about cohomology operations is argued in \cite{HoyoisKO17}*{\S 2.4} and \cite{Hoyois15}*{\S 5.2}.
\end{remark}

\begin{corollary}\label{cor:ReduceBase}
Let $f \colon S \to T$ be a map of base schemes such that the structure map $f^* \MM\F_{\ell}^{T} \ral{\cong} \MM\F_{\ell}^S$ in $\SH(S)$ is an isomorphism. If the map $\Psi^T$ in $\hMod{\MM\F_{\ell}^T}$ is an isomorphism (resp. split monomorphism), then so is the map $\Psi^S$ in $\hMod{\MM\F_{\ell}^S}$. 
\end{corollary}

A similar base change argument was used in the proof of \cite{HoyoisKO17}*{Theorem~5.6}.

\begin{lemma}\label{lem:RingCoefficients}
Let $f \colon S \to T$ be a map of base schemes such that the structure map $f^* \MM\Z_T \ral{\cong} \MM\Z_S$ in $\SH(S)$ is an isomorphism. Then for any abelian group $A$, the structure map $f^* \MM A_T \ral{\cong} \MM A_S$ in $\SH(S)$ is an isomorphism.
\end{lemma}

\begin{proof}
See \cite{Hoyois15}*{Proposition~4.13~(2)} and the discussion thereafter. 
For example, $M(\Z/n)_S$ is given by the cofiber $(\MM\Z_S) / n$.
\end{proof}

\begin{lemma}\label{lem:ExtnPerfect}
\cite{Hoyois15}*{Lemma~A.2}
Let $\kk$ be a perfect field and $\kk \to \mathbb{L}$ a field extension. Then the induced map of schemes $\Spec(\mathbb{L}) \to \Spec(\kk)$ is essentially smooth.
\end{lemma}

\begin{corollary}\label{cor:ReduceFp}
If $\Psi^{\F_p}$ is an isomorphism, then $\Psi^{\kk}$ is an isomorphism for every field $\kk$ of characteristic $p$.
\end{corollary}

\begin{proof}
By Lemma~\ref{lem:ExtnPerfect}, the map of schemes $f \colon \Spec(\kk) \to \Spec(\F_p)$ is essentially smooth. By Theorem~\ref{thm:PullbackHZ} \eqref{item:PullbackSmooth}, $f^* \MM\Z_{\F_p} \ral{\cong} \MM\Z_{\kk}$ is an isomorphism. By Lemma~\ref{lem:RingCoefficients}, $f^* \MM\F_{\ell}^{\F_p} \ral{\cong} \MM\F_{\ell}^{\kk}$ is also an isomorphism. The statement now follows from Corollary~\ref{cor:ReduceBase}.
\end{proof}

\section{Discrete valuation ring of mixed characteristic}\label{sec:DVR}

\begin{lemma}
Let $\kk$ be a field of characteristic $p$. Then there exists a discrete valuation ring $D$ having $\kk$ as residue field, and a fraction field $\KK = \Frac(D)$ of characteristic zero. 
\end{lemma}

\begin{proof}
If $\kk$ is a perfect field, then the $p$-adic Witt ring $D = W_{p^{\infty}}(\kk)$ satisfies the desired properties \cite{Hazewinkel09}*{Theorem~6.19, 6.20}. 
More generally, whether $\kk$ is perfect or not, the Cohen ring of $\kk$ satisfies the desired properties \cite{Hazewinkel09}*{6.23}.
\end{proof}

\begin{example}
In the case $\kk = \F_p$, one can take the $p$-adic integers $D = W_{p^{\infty}}(\F_p) = \Z_p$, with fraction field the $p$-adic rationals $\KK = \Q_p$.
\end{example}

\begin{notation}\label{nota:RingMaps}
Consider the ring maps and induced maps of affine schemes:
\[
\xymatrix{
& \KK & & & & \Spec(\KK) \ar@{^{(}->}[dl]_-{j}^-{\text{open}} \\
D \ar[ur] \ar@{->>}[dr] & & \ar@{~>}[r] & & \Spec(D) & \\
& \kk & & & & \Spec(\kk). \ar@{^{(}->}[ul]^-{i}_-{\text{closed}} \\
}
\]
\end{notation}

That is, $i \colon \Spec(\kk) \inj \Spec(D)$ is the inclusion of the closed point, while $j \colon \Spec(\KK) \inj \Spec(D)$ is its open complement, the generic point. The strategy of the proof is to go through $\Spec(D)$ to transfer a result which is known over $\Spec(\KK)$ to the desired result over $\Spec(\kk)$.

\begin{notation}
Let $\hatt{M}\Z$ denote the version of the Eilenberg--MacLane spectrum constructed in \cite{Spitzweck18}. 
There is a variant $\hatt{M}A$ for any abelian group $A$. %
\end{notation}

Recall some salient facts about $\hatt{M}\Z$ proved in \cite{Spitzweck18}. 

\begin{theorem}\label{thm:HZ-Spitzweck}
\begin{enumerate}
\item $\hatt{M}\Z$ is a motivic $E_{\infty}$ ring spectrum \cite{Spitzweck18}*{\S 4.4}.
\item The family $\hatt{M}\Z = \{ \hatt{M}\Z_S \}_S$ is Cartesian, i.e., all structure maps $f^* \hatt{M}\Z_T \ral{\cong} \hatt{M}\Z_S$ are isomorphisms. \cite{Spitzweck18}*{Theorem~8.25}.
\item Over a field $\kk$, $\hatt{M}\Z$ agrees with the usual motivic Eilenberg--MacLane spectrum, i.e., there is a canonical isomorphism $\hatt{M}\Z_{\kk} \cong \MM\Z_{\kk}$ in $\Ho (E_{\infty}-\Spt(\kk))$ \cite{Spitzweck18}*{Theorem~8.22}. 
\item The $\base$-fibered category 
$\hMod{\hatt{M}\Z}$ (as in Definition~\ref{def:DM}) satisfies the six functor formalism \cite{Spitzweck18}*{\S 9}. 
\end{enumerate}
\end{theorem}

\begin{corollary}\label{cor:EMpoint}
Using Notation~\ref{nota:RingMaps}, there are canonical isomorphisms %
\begin{enumerate}
\item $j^* \hatt{M}\F_{\ell}^{D} \cong \MM\F_{\ell}^{\KK}$ in $\Ho (E_{\infty}-\Spt(\KK))$
\item $i^* \hatt{M}\F_{\ell}^D \cong \MM\F_{\ell}^{\kk}$ in $\Ho (E_{\infty}-\Spt(\kk))$.
\end{enumerate}
\end{corollary}

\begin{warning}\label{war:HZ-Spitzweck}
In view of Corollary~\ref{cor:EMpoint}, from now on we will usually drop the hats from the notation $\hatt{M}\Z$, though all the arguments will be using $\hatt{M}\Z$ instead of Voevodsky's Eilenberg--MacLane spectrum $\MM\Z$, as in Notation~\ref{nota:EMspectra}. 
Over the base ring $D$, the spectrum $\hatt{M}\Z_D$ is \emph{not} known to be equivalent to $\MM\Z_D$, though conjectured to be.
\end{warning}

\begin{lemma}\label{lem:MotivicH21}
For any smooth scheme $S$ over $\kk$, we have the following motivic cohomology groups
\[
\begin{cases}
H^{2,1} (S; \Z) = \Pic(S) \\
H^{2,1} (S; \Z/\ell) = \Pic(S) \ot_{\Z} \Z/\ell. \\
\end{cases}
\]
\end{lemma}

\begin{proof}
The first part is given in \cite{MazzaVW06}*{Corollary~4.2}. The second part follows from the universal coefficient exact sequence
\begin{equation}\label{eq:UnivCoefSeq}
\xymatrix{
0 \ar[r] & H^{2,1}(S;\Z) \ot_{\Z} \Z/\ell \ar[r] & H^{2,1}(S;\Z/\ell) \ar[r] & \Tor_1^{\Z} \left( H^{3,1}(S;\Z), \Z/\ell \right) \ar[r] & 0 \\
}
\end{equation}
and the equality $H^{3,1} (S; \Z) = 0$.
\end{proof}

In what follows, we will use some facts from \cite{Ayoub07I}, %
and the fact that $\hMod{\hatt{M}\F_{\ell}}$ satisfies the six functor formalism (Theorem~\ref{thm:HZ-Spitzweck}). 
See also %
\cite{CisinskiD19}*{\S 2.3} 
for more details on the localization property. 
The new key input is the isomorphism of $\MM\Z_{\kk}$-modules $i^! \MM\Z_{D} \cong \Si^{-2,-1} \MM\Z_{\kk}$, which will be proved in Corollary~\ref{cor:ShriekHZ}.

\begin{lemma}\label{lem:EMopenClosed}
Let $\ell$ be any prime number. Then there is a splitting
\[
i^* j_* \MM\F_{\ell}^{\KK} \simeq \MM\F_{\ell}^{\kk} \op \Si^{-1,-1} \MM\F_{\ell}^{\kk} 
\]
in $\hMod{\MM\F_{\ell}^{\kk}}$.%
\end{lemma}

\begin{proof}
For any object $E$ in $\hMod{\MM\F_{\ell}^D}$, %
by the localization property, there is an exact triangle
\begin{equation}\label{eq:ExactTriD}
\xymatrix{
i_! i^! E \ar[r]^-{\ep} & E \ar[r]^-{\eta} & j_* j^* E \\
}
\end{equation}
in $\hMod{\MM\F_{\ell}^D}$, natural in $E$. %
Since $i$ is a closed immersion, we have $i_! = i_*$. 
Applying the exact functor $i^*$ yields an exact triangle
\[
\xymatrix{
i^* i_* i^! E \ar[r]^-{i^* \ep} & i^* E \ar[r]^-{i^* \eta} & i^* j_* j^* E \\
}
\]
in $\hMod{\MM\F_{\ell}^{\kk}}$. %
Since $i$ is a closed immersion, the counit $i^* i_* \to 1$ is an isomorphism (by the localization property), and we can rewrite the exact triangle as 
\begin{equation}\label{eq:ExactTriK}
\xymatrix{
i^! E \ar[r] & i^* E \ar[r] & i^* j_* j^* E. \\
}
\end{equation}
Now take $E = \MM\F_{\ell}^{D}$. Using the isomorphisms $j^* \MM\F_{\ell}^{D} = \MM\F_{\ell}^{\KK}$ and $i^* \MM\F_{\ell}^{D} = \MM\F_{\ell}^{\kk}$ of Corollary~\ref{cor:EMpoint}, the exact triangle~\eqref{eq:ExactTriK} in $\hMod{\MM\F_{\ell}^{\kk}}$ %
becomes 
\begin{equation}\label{eq:ExactTriHFp}
\xymatrix{
i^! \MM\F_{\ell}^{D} \ar[r] & \MM\F_{\ell}^{\kk} \ar[r] & i^* j_* \MM\F_{\ell}^{\KK}. \\
}
\end{equation}
Using the isomorphism $i^! \MM\F_{\ell}^{D} \cong \Si^{-2,-1} \MM\F_{\ell}^{\kk}$ from Corollary~\ref{cor:ShriekHZ}, the exact triangle~\eqref{eq:ExactTriHFp} becomes
\[
\xymatrix{
\Si^{-2,-1} \MM\F_{\ell}^{\kk} \ar[r] & \MM\F_{\ell}^{\kk} \ar[r] & i^* j_* \MM\F_{\ell}^{\KK}. \\
}
\]
Next, we claim that the map $\Si^{-2,-1} \MM\F_{\ell}^{\kk} \to \MM\F_{\ell}^{\kk}$ must be zero. Indeed, such $\MM\F_\ell$-module maps are given by
\begin{align*}
\hMod{\MM\F_{\ell}^{\kk}} \left( \Si^{-2,-1} \MM\F_{\ell}^{\kk}, \MM\F_{\ell}^{\kk} \right) &= \hMod{\MM\F_{\ell}^{\kk}} \left( \MM\F_{\ell}^{\kk}, \Si^{2,1} \MM\F_{\ell}^{\kk} \right) \\
&= \SH(\kk) \left( \unit_{\kk}, \Si^{2,1} \MM\F_{\ell}^{\kk} \right) \\
&= H^{2,1} (\Spec(\kk); \F_{\ell}) \\
&= 0.
\end{align*}
Here we used the fact that motivic Eilenberg--MacLane spectra represent motivic cohomology \cite{Voevodsky98}*{\S 6.2} \cite{Hoyois15}*{Theorem~4.24}, along with Lemma~\ref{lem:MotivicH21} and the equality $\Pic \left( \Spec(\kk) \right) = \Pic(\kk) = 0$. Thus, the exact triangle~\eqref{eq:ExactTriHFp} splits, which yields the isomorphism
\begin{align*}
i^* j_* \MM\F_{\ell}^{\KK} &\simeq \MM\F_{\ell}^{\kk} \op \Si^{1,0} \Si^{-2,-1} \MM\F_{\ell}^{\kk} \\
&= \MM\F_{\ell}^{\kk} \op \Si^{-1,-1} \MM\F_{\ell}^{\kk}
\end{align*}
as claimed.
\end{proof}

\section{Retract of the dual Steenrod algebra}\label{sec:RetractSteenrod}

In this section, we prove our first main result.

\begin{theorem}\label{thm:RetractSteenrod}
Assume that the base scheme $S$ is essentially smooth over a field $\kk$ of characteristic $p > 0$. 
Then the map $\Psi^{S}$ in $\hMod{\MM\F_p^S}$ displayed in Equation~\eqref{eq:Decompo} is a split monomorphism, i.e., admits a retraction. 
\end{theorem}

\begin{proof}
By Theorem~\ref{thm:PullbackHZ}, Lemma~\ref{lem:RingCoefficients}, and Corollary~\ref{cor:ReduceBase}, we may assume $S = \kk$. 
Consider the diagram in $\hMod{\MM\F_p^{\kk}}$
\begin{equation}\label{eq:DiagRetract}
\xymatrix @R-0.4pc @C=\bigcol {
\bigoplus_{\al} \Si^{p_{\seq}, q_{\seq}} \MM\F_p^{\kk} \ar[r]^-{\Psi^{\kk}} & \MM\F_p^{\kk} \sm_{\kk} \MM\F_p^{\kk} \\
\bigoplus_{\seq} \Si^{p_{\seq}, q_{\seq}} i^* \MM\F_p^{D} \ar[u]^{\cong} \ar[d]_{\cong} & i^* \MM\F_p^{D} \sm_{\kk} i^* \MM\F_p^{D} \ar[u]_{\cong} \ar[d]^{\cong} \\
i^* ( \bigoplus_{\seq} \Si^{p_{\seq}, q_{\seq}} \MM\F_p^{D} ) \ar[d]_{i^* \eta} \ar[r]^-{i^* \Psi^{D}} & i^* ( \MM\F_p^{D} \sm_{D} \MM\F_p^{D} ) \ar[d]^{i^* \eta} \\
i^* j_* j^* ( \bigoplus_{\seq} \Si^{p_{\seq}, q_{\seq}} \MM\F_p^{D} ) \ar[r]^-{i^* j_* j^* \Psi^{D}} & i^* j_* j^* ( \MM\F_p^{D} \sm_{D} \MM\F_p^{D} ) \\
i^* j_* ( \bigoplus_{\seq} \Si^{p_{\seq}, q_{\seq}} j^* \MM\F_p^{D} ) \ar[d]_{\cong} \ar[u]^{\cong} & i^* j_* ( j^* \MM\F_p^{D} \sm_{\KK} j^* \MM\F_p^{D} ) \ar[u]_{\cong} \ar[d]^{\cong} \\
i^* j_* ( \bigoplus_{\seq} \Si^{p_{\seq}, q_{\seq}} \MM\F_p^{\KK} ) \ar[r]^-{i^* j_* \Psi^{\KK}}_{\cong} & i^* j_* ( \MM\F_p^{\KK} \sm_{\KK} \MM\F_p^{\KK} ) \\
\bigoplus_{\seq} \Si^{p_{\seq}, q_{\seq}} i^* j_* \MM\F_p^{\KK}. \ar[u]^{\cong} & \\
}
\end{equation}
The top and bottom rectangles commute, by the identifications $i^* \Psi^{D} = \Psi^{\kk}$ and $j^* \Psi^{D} = \Psi^{\KK}$ from Lemma~\ref{lem:CompatDecompo}. The middle rectangle commutes by naturality of the unit map $\eta$. Define the map
\[
\phy' \colon \MM\F_p^{\kk} \sm_{\kk} \MM\F_p^{\kk} \to \bigoplus_{\seq} \Si^{p_{\seq}, q_{\seq}} i^* j_* \MM\F_p^{\KK}
\]
as the composite starting from the top right corner, going all the way down, then left, then down again. We claim that the composite
\[
\phy' \circ \Psi^{\kk} \colon \bigoplus_{\seq} \Si^{p_{\seq}, q_{\seq}} \MM\F_p^{\kk} \to \bigoplus_{\seq} \Si^{p_{\seq}, q_{\seq}} i^* j_* \MM\F_p^{\KK}
\]
is equal to the map $\bigoplus_{\seq} \Si^{p_{\seq}, q_{\seq}} i^* \eta$, for the map $i^* \eta \colon \MM\F_p^{\kk} \to i^* j_* \MM\F_p^{\KK}$. This follows from the commutative diagram in $\hMod{\MM\F_p^D}$
\[
\xymatrix @C+1pc {
\bigoplus_{\seq} \Si^{p_{\seq}, q_{\seq}} \MM\F_p^{D} \ar[d]_-{\bigoplus_{\seq} \Si^{p_{\seq}, q_{\seq}} \eta} \ar@{=}[r] & \bigoplus_{\seq} \Si^{p_{\seq}, q_{\seq}} \MM\F_p^{D} \ar[d]^-{\eta} \\
\bigoplus_{\seq} \Si^{p_{\seq}, q_{\seq}} j_* j^* \MM\F_p^{D} \ar[dd]_{\cong} \ar[dr]^-{\cong} \ar[r]^-{\cong} & j_* j^* ( \bigoplus_{\seq} \Si^{p_{\seq}, q_{\seq}} \MM\F_p^{D} ) \\
 & j_* ( \bigoplus_{\seq} \Si^{p_{\seq}, q_{\seq}} j^* \MM\F_p^{D} ) \ar[d]^{\cong} \ar[u]_{\cong} \\
\bigoplus_{\seq} \Si^{p_{\seq}, q_{\seq}} j_* \MM\F_p^{\KK} \ar[r]^-{\cong} & j_* ( \bigoplus_{\seq} \Si^{p_{\seq}, q_{\seq}} \MM\F_p^{\KK} ). \\
}
\]
By Lemma~\ref{lem:EMopenClosed}, the map $i^* \eta \colon \MM\F_p^{\kk} \to i^* j_* \MM\F_p^{\KK}$ in $\hMod{\MM\F_p^{\kk}}$ admits a retraction, say, $r \colon i^* j_* \MM\F_p^{\KK} \to \MM\F_p^{\KK}$. Define $\phy$ as the composite
\[
\xymatrix @C=\bigcol {
\MM\F_p^{\kk} \sm_{\kk} \MM\F_p^{\kk} \ar[r]^-{\phy'} & \bigoplus_{\seq} \Si^{p_{\seq}, q_{\seq}} i^* j_* \MM\F_p^{\KK} \ar[r]^-{\bigoplus_{\seq} \Si^{p_{\seq}, q_{\seq}} r} & \bigoplus_{\seq} \Si^{p_{\seq}, q_{\seq}} \MM\F_p^{\kk}. \\
}
\]
Then $\phy$ satisfies
\begin{align*}
\phy \circ \Psi^{\kk} &= (\bigoplus_{\seq} \Si^{p_{\seq}, q_{\seq}} r) \circ \phy' \circ \Psi^{\kk} \\
&= \bigoplus_{\seq} \Si^{p_{\seq}, q_{\seq}} (r \circ i^* \eta) \\
&= \id,
\end{align*}
providing the desired retraction.
\end{proof}

\section{Retract of the slices of \texorpdfstring{$\MGL$}{MGL}}\label{sec:RetractSlices}

\subsection{The slice filtration}

Background on the slice filtration can be found in \cite{Voevodsky02ope}*{\S 2}, \cite{Spitzweck12}*{\S 3}, \cite{GutierrezRSO12}*{\S 3}, \cite{RondigsO16}*{\S 2}, and \cite{Bachmann17}*{\S 1}. We fix some notation and terminology.

Let $i_n \colon \Si^{2n,n} \SH(S)^{\eff} \inj \SH(S)$ denote the inclusion of the subcategory of $n$-effective spectra, i.e., the localizing triangulated subcategory generated by objects $\Si^{2n,n} \Si^{\infty}_T X_+$ for $X \in \Sm_S$. Let $r_n \colon \SH(S) \to \Si^{2n,n} \SH(S)^{\eff}$ denote the right adjoint to the inclusion $i_n$, and $f_n \dfn i_n r_n$ the associated comonad, with counit $\cov_n \colon f_n E \to E$, that is, the $n$-effective cover functor. These form a tower over $E$
\[
\ldots \to f_n E \to f_{n-1} E \to \ldots \to E.
\]
The $n$\textsuperscript{th} slice sits in a natural distinguished triangle
\begin{equation}\label{eq:SliceTriangle}
\xymatrix{
f_{n+1} E \ar[r] & f_n E \ar[r] & s_n E \ar[r] & \Si^{1,0} f_{n+1} E.
}
\end{equation}
Let $\tr^n \colon E \to f^n E$ denote the natural truncation map, %
sitting in 
the natural distinguished triangle
\begin{equation}\label{eq:TruncationTriangle}
\xymatrix{
f_{n+1} E \ar[r]^-{\cov_n} & E \ar[r]^-{\tr^n} & f^{n} E \ar[r] & \Si^{1,0} f_{n+1} E. \\
}
\end{equation}
These truncations form a tower under $E$
\[
E \to \ldots \to f^n E \to f^{n-1} E \to \ldots.
\]
Successive truncations sit in a natural distinguished triangle
\begin{equation}\label{eq:SliceTruncation}
\xymatrix{
s_{n} E \ar[r] & f^{n} E \ar[r] & f^{n-1} E \ar[r] & \Si^{1,0} s_{n} E.
}
\end{equation}
Applying the exact functor $s_i$ to the exact triangle~\eqref{eq:TruncationTriangle} yields a description of the slices of the $n$-truncation $f^n E$. For $i > n$, we have
\[
s_i (f^n E) = 0
\]
since the map $s_i (\cov_n)$ is an isomorphism. For $i \leq n$, the truncation map $\tr^n \colon E \to f^n E$ induces isomorphisms on slices
\[
s_i E \ral{\cong} s_i (f^n E)
\]
since $s_i (f^{n+1} E) = 0$ holds. An object $E$ is called \Def{$n$-truncated} if the truncation map $\tr^n \colon E \to f^n E$ is an isomorphism.

The following can be found in \cite{Voevodsky02ope}*{proof of Lemma~4.2}, %
\cite{RondigsSO18}*{Lemma~3.1}, and \cite{Hoyois15}*{\S 8.3}.

\begin{lemma}\label{lem:Exhaustive}
For every object $E$ if $\SH(S)$, the slice filtration is exhaustive, i.e., the natural map
\[
\hocolim_{n \to -\infty} f_n E \ral{\cong} E
\]
is an isomorphism.
\end{lemma}

\begin{lemma}\label{lem:SliceShift}
For all $E \in \SH(S)$ and $n \in \Z$, there are natural isomorphisms in $\SH(S)$ as follows.
\begin{enumerate}
\item
\[
\begin{cases}
f_n(\Si^{0,1} E) \cong \Si^{0,1}(f_{n-1} E) \\
f_n(\Si^{1,0} E) \cong \Si^{1,0}(f_n E) \\
\end{cases}
\]
and thus $f_n(\Si^{p,q} E) \cong \Si^{p,q}(f_{n-q} E)$ for all $p,q \in \Z$.
\item 
\[
\begin{cases}
s_n(\Si^{0,1} E) \cong \Si^{0,1}(s_{n-1} E) \\
s_n(\Si^{1,0} E) \cong \Si^{1,0}(s_n E) \\
\end{cases}
\]
and thus $s_n(\Si^{p,q} E) \cong \Si^{p,q}(s_{n-q} E)$ for all $p,q \in \Z$.
\end{enumerate}
\end{lemma}

\begin{proof}
This follows from the same argument as in \cite{RondigsO16}*{Lemma~2.1}.
\end{proof}

For $E \in \SH(S)$ and $p,q \in \Z$, recall that the \Def{homotopy sheaf} $\ul{\pi}_{p,q}(E)$ is defined as the Nisnevich sheafification of the presheaf of abelian groups $\ul{\pi}^{\pre}_{p,q}(E)$ given by
\[
X \mapsto [\Si^{p,q} \Si^{\infty} X_+, E]
\]
for $X \in \Sm_S$ \cite{Spitzweck20}*{\S 4}.

\begin{lemma}\label{lem:HomotGpsCover}
Let $E \in \SH(S)$ and $n \in \Z$.
\begin{enumerate}
\item The $n$-effective cover map $\cov_n \colon f_n E \to E$ induces an isomorphism on homotopy %
sheaves
\[
\ul{\pi}_{p,q} (f_n E) \ral{\cong} \ul{\pi}_{p,q} E
\]
for $q \geq n$.
\item The $n$-truncation $f^n E$ satisfies
\[
\ul{\pi}_{p,q} (f^n E) = 0
\]
for $q \geq n+1$.
\end{enumerate}
\end{lemma}

\begin{proof}
For part~(1), the sphere $S^{p,q} = \Si^{p,q} \unit_S$ is $q$-effective, in particular $n$-effective assuming $q \geq n$. For any $X \in \Sm_S$, the spectrum $\Si^{p,q} \Si^{\infty} X_+$ is also $q$-effective. The universal property of the counit map $\cov_n$ yields the %
isomorphism
\[
[\Si^{p,q} \Si^{\infty} X_+, f_n E] \ral{\cong} [\Si^{p,q} \Si^{\infty} X_+, E],
\]
so that $\cov_n$ induces an isomorphism of presheaves $\ul{\pi}^{\pre}_{p,q}(f_n E) \ral{\cong} \ul{\pi}^{\pre}_{p,q}(E)$. Sheafifying then yields the claimed isomorphism.

Part~(2) follows from part~(1), the exact triangle~\eqref{eq:TruncationTriangle}, and the long exact sequence of homotopy sheaves.
\end{proof}

\begin{lemma}\label{lem:CoverTrunc}
There are natural isomorphisms $s_n \cong f_n f^n \cong f^n f_n$ of functors $\SH(S) \to \SH(S)$, which moreover identify the commutative diagram
\[
\xymatrix{
f_n E \ar[d]_{\cov_n} \ar[r]^-{f_n(\tr^n)} & f_n f^n E \ar[d]^{\cov_n} \\
E \ar[r]^-{\tr^n} & f^n E \\
}
\]
with
\[
\xymatrix{
f_n E \ar[d]_{\cov_n} \ar[r]^-{\tr^n} & f^n f_n E \ar[d]^{f^n (\cov_n)} \\
E \ar[r]^-{\tr^n} & f^n E. \\
}
\]
\end{lemma}

\begin{proof}
See \cite{GutierrezRSO12}*{\S 2.3}.
\end{proof}

\begin{lemma}\label{lem:SliceSlice}
For every $E \in \SH(S)$ and $m,n \in \Z$, there are natural isomorphisms
\[
s_m s_n E = \begin{cases}
s_n E &\text{if } m=n \\
0 &\text{if } m \neq n.
\end{cases}
\]
\end{lemma}

The commutative diagram from Lemma~\ref{lem:CoverTrunc}
\[
\xymatrix{
f_n E \ar[d]_{\cov_n} \ar[r]^-{f_n(\tr^n)} & s_n E \ar[d]^{\cov_n} \\
E \ar[r]^-{\tr^n} & f^n E \\
}
\]
induces on homotopy groups $\pi_{p,q}$ with $q \geq n$ the diagram
\begin{equation}\label{eq:DiagHomotGps}
\xymatrix{
\pi_{p,q} (f_n E) \ar[d]_{\cong} \ar[r]^-{} & \pi_{p,q} (s_n E) \ar[d]^{\cong} \\
\pi_{p,q} E \ar[r]^-{} & \pi_{p,q} (f^n E), \\
}
\end{equation}
using the isomorphism from Lemma~\ref{lem:HomotGpsCover}. 

\begin{notation}\label{nota:HomotSlice}
For $q \geq n$, let 
\[
\si \colon \pi_{p,q} E \to \pi_{p,q} (s_n E)
\]
denote the natural map of abelian groups defined by the diagram~\eqref{eq:DiagHomotGps}. This was used for instance in \cite{Voevodsky02ope}*{\S 3.2}.
\end{notation}

\begin{lemma}\label{lem:PullbackSlice}
Let $g \colon S \to T$ be a map of base schemes and let $n \in \Z$.
\begin{enumerate}
\item If $E \in \SH(T)$ is $n$-effective, then $g^*E \in \SH(S)$ is $n$-effective. In particular, there is a canonical natural transformation $g^* f_n E \to f_n g^* E$ which in turn induces a natural transformation $g^* s_n E \to s_n g^* E$.
\item If $g$ is essentially smooth, then the canonical transformations are isomorphisms $g^* f_n \cong f_n g^*$ and $g^* s_n \cong s_n g^*$.  
\end{enumerate}
\end{lemma}

\begin{proof}
See \cite{Hoyois15}*{Remark~4.20} or \cite{RondigsO16}*{Theorem~2.5}.
\end{proof}

\subsection{Some facts about slices}

We recall some facts about known slices. More details and references can be found in \cite{GutierrezRSO12}*{\S 6.1}.

\begin{theorem}\label{thm:SlicesHZ}
Assume that the base scheme $S$ is essentially smooth over a field.
\begin{enumerate}
\item The Eilenberg--MacLane spectrum $\MM\Z_S$ has slices $s_0 \MM\Z = \MM\Z$ and $s_n \MM\Z = 0$ for $n \neq 0$. (Voevodsky's Conjecture~1.)
\item The unit map $u \colon \unit \to \MM\Z$ induces an isomorphism on zeroth slices $s_0 \unit \ral{\cong} s_0 \MM\Z \cong \MM\Z$. (Voevodsky's Conjecture~10.)
\item The algebraic $K$-theory spectrum $\KGL_S$ has slices $s_n \KGL = \Si^{2n,n} \MM\Z$ for $n \in \Z$. (Voevodsky's Conjecture~7.)
\end{enumerate}
\end{theorem}

The isomorphism $s_0 \unit \cong \MM\Z$ was proved in \cite{Voevodsky04} for fields of characteristic zero, generalized to perfect fields in \cite{Levine08}, 
from which one obtains the stated generality \cite{Hoyois15}*{Remark~4.20}. 
The result on the slices of $\KGL$ is also due to Voevodsky and Levine; see \cite{RondigsO16}*{Theorem~4.1}.

The slice filtration is compatible with multiplicative structures, in the following sense.

\begin{theorem}\label{thm:SliceMultiplicative}
Let $E$ be a motivic $E_{\infty}$ ring spectrum and $M$ an $E$-module spectrum. Then for every $q \in \Z$, the cover $f_q M$ is canonically an $f_0 E$-module and the slice $s_q E$ is canonically an $s_0 E$-module \cite{Pelaez11} \cite{GutierrezRSO12}*{Corollary~5.18, \S 6(v)}. 
\end{theorem}

In particular, for every motivic spectrum $M$ over $S$, the slice $s_q M$ is canonically an $\MM\Z$-module, assuming the isomorphism $s_0 \unit \cong \MM\Z$.

\begin{lemma}\label{lem:SelfMapsHZ}
Assume that the base scheme $S$ is connected and essentially smooth over a field. Then the graded self maps of $\MM\Z$ in $\SH(S)$ are given by
\[
[\MM\Z, \Si^{s,0} \MM\Z] = \begin{cases}
\Z &\text{if } s = 0 \\
0 &\text{otherwise.} \\
\end{cases} 
\] 
\end{lemma}

See also \cite{Voevodsky02ope}*{Conjecture~12} and the remark thereafter.

\begin{proof}
By Theorem~\ref{thm:SlicesHZ}, the unit map $u \colon \unit \to \MM\Z$ induces an isomorphism on zeroth slices $s_0 \unit \ral{\cong} s_0 \MM\Z \cong \MM\Z$. Since the sphere spectrum $\unit = f_0 \unit$ is effective, the unit map can be identified with the $0$-truncation $\tr^0 \colon f_0 \unit \to s_0 \unit \cong \MM\Z$. Since $\Si^{s,0} \MM\Z = f^0 (\Si^{s,0} \MM\Z)$ is $0$-truncated, precomposition by the unit map $\unit \to \MM\Z$ induces the isomorphism
\[
[\MM\Z, \Si^{s,0} \MM\Z] \ral{\cong} [\unit, \Si^{s,0} \MM\Z] = \pi_{-s,-0} \MM\Z.
\] 
Since $\MM\Z$ represents motivic cohomology \cite{Hoyois15}*{Theorem~4.24}, 
we further obtain the isomorphism
\[
[\unit, \Si^{s,0} \MM\Z] = [\Si^{\infty}_+ S, \Si^{s,0} \MM\Z] = H^{s,0}(S;\Z).
\]
Since $S$ is %
smooth over a field, we have the motivic cohomology groups
\[
H^{s,0}(S;\Z) = \begin{cases}
\Z^{\Conn(S)} &\text{if } s = 0 \\
0 &\text{otherwise,}
\end{cases}
\]
by \cite{MazzaVW06}*{Corollary~4.2}.
\end{proof}

That is, $H^{0,0}(S;\Z) = \prod_{\Conn(S)} \Z$ is the ring of $\Z$-valued functions on the set $\Conn(S)$ of connected components of the scheme $S$. 

\begin{corollary}\label{cor:EndoHZ}
Over a base field, taking the homotopy group $\pi_{0,0}$ induces a ring isomorphism
\[
\xymatrix{
[\MM\Z, \MM\Z] \ar[r]^-{\pi_{0,0}}_-{\cong} & \Hom_{\Z}(\Z,\Z). \\
}
\]
Here, we use the %
isomorphism $\Z \ral{\cong} \pi_{0,0} \MM\Z$ induced by the unit map $\unit \to \MM\Z$.

In particular, a map $\MM\Z \to \MM\Z$ in $\SH(S)$ is an isomorphism if and only if it induces a surjection on $\pi_{0,0}$.
\end{corollary}

\subsection{Slices of $\MGL$}

Let
\[
L_* \cong \Z[x_1, x_2, \ldots], \quad \abs{x_i} = i
\]
denote the Lazard ring. Here the grading is half of the grading that appears in topology: 
$L_* \cong \MU_{2*}$. %
Recall that there is a canonical ring homomorphism
\[
\la_S \colon L_{*} \to \pi_{2*,*} \MGL_S,
\]
described for instance in \cite{Voevodsky02ope}*{\S 3.2}, \cite{Spitzweck10}*{\S 4}, and \cite{Hoyois15}*{\S 6.1}. 
For fixed $n \geq 0$, the map of abelian groups $\la_{S} \colon L_{n} \to \pi_{2n,n} \MGL_S$ corresponds to a map in $\SH(S)$
\[
\tild{\la}_S \colon \Si^{2n,n} \unit_S \ot_{\Z} L_{n} \to \MGL_S.
\]
This ring homomorphism $\la$ is compatible with change of base scheme in the following sense.

\begin{lemma}\label{lem:LazardMap}
Let $g \colon S \to T$ be a map of base schemes.
\begin{enumerate}
\item The natural map $g^* \MGL_T \ral{\cong} \MGL_S$ in $\SH(S)$ is an isomorphism.
\item  The following diagram in $\SH(S)$ commutes:
\[
\xymatrix{
g^* ( \Si^{2n,n} \unit_T \ot_{\Z} L_{n} ) \ar[d]_{\cong} \ar[r]^-{g^* \tild{\la}_T} & g^* \MGL_T \ar[d]^{\cong} \\
\Si^{2n,n} \unit_S \ot_{\Z} L_{n} \ar[r]^-{\tild{\la}_S} & \MGL_S. \\
}
\]
\end{enumerate}
\end{lemma}

\begin{proof}
Part~(1) follows from \cite{NaumannSO09}*{Proposition~8.5} with $M_* = \MU_*$. 
See also \cite{Deglise18ori}*{Example~1.2.3(3)} and the references therein. 
Part~(2) follows from the construction described in \cite{Spitzweck12}*{\S 5}. See also \cite{Hoyois15}*{\S 6.1, after Definition~7.1}.
\end{proof}

The map of abelian groups
\[
\xymatrix{
L_{n} \ar[r]^-{\la_S} & \pi_{2n,n} \MGL_S \ar[r]^-{\si} & \pi_{2n,n} (s_n \MGL_S) \\
}
\]
corresponds to a map in $\SH(S)$
\begin{equation}\label{eq:ApproxSliceSH}
\Si^{2n,n} \unit_S \ot_{\Z} L_{n} \to s_n \MGL_S.
\end{equation}

\begin{notation}\label{nota:ApproxSliceDM}
Assume that the base scheme $S$ satisfies $s_0 \unit_S \cong \MM\Z_S$, as in Theorem~\ref{thm:SlicesHZ}. Since the slice $s_n \MGL_S$ is canonically an $\MM\Z_S$-module (by Theorem~\ref{thm:SliceMultiplicative}), the map~\eqref{eq:ApproxSliceSH} in $\SH(S)$ corresponds by adjunction to a map in $\hMod{\MM\Z_{S}}$ %
\begin{equation}\label{eq:ApproxSliceDM}
\xymatrix{
\Si^{2n,n} \MM\Z_{S} \ot_{\Z} L_{n} \ar[r]^-{\psi_S} & s_n \MGL_{S}. \\
}
\end{equation}
\end{notation}

\begin{conjecture}\label{conj:SliceMGL}
The map %
$\psi_S \colon \Si^{2n,n} \MM\Z_{S} \ot_{\Z} L_{n} \to s_n \MGL_{S}$ 
in $\hMod{\MM\Z_S}$ is an isomorphism \cite{Voevodsky02ope}*{Conjecture~5}.
\end{conjecture}

By construction, the map $\psi_S$ is compatible with the complex orientation of $\MGL$, i.e., it satisfies the compatibility condition described after~\cite{Voevodsky02ope}*{Conjecture~5}. More precisely, the following diagram of abelian groups commutes:
\begin{equation}\label{eq:CompatLazard}
\xymatrix{
L_{n} \ar[d]_-{\la} \ar[r]^-{\text{unit}} & \pi_{2n,n} (\Si^{2n,n} \MM\Z \ot_{\Z} L_{n}) \ar[d]^-{\pi_{2n,n} \psi} \\
\pi_{2n,n} \MGL \ar[r]^-{\si} & \pi_{2n,n} (s_n \MGL). \\
}
\end{equation}
The map at the top is induced by the natural map $\Z \to \pi_{0,0} \MM\Z_S$, cf.\ Lemma~\ref{lem:SelfMapsHZ}. 
By adjunction, and by construction of the natural map $\si \colon \pi_{2n,n} E \to \pi_{2n,n} (s_n E)$ in Notation~\ref{nota:HomotSlice}, commutativity of the diagram~\eqref{eq:CompatLazard} is equivalent to commutativity of the following diagram in $\SH(S)$:
\begin{equation}\label{eq:CompatLazardSH}
\xymatrix @R-0.4pc {
\Si^{2n,n} \unit \ot_{\Z} L_{n} \ar[dd]_-{\tild{\la}} \ar[r]^-{u} & \Si^{2n,n} \MM\Z \ot_{\Z} L_{n} \ar[d]^-{\psi} \\
& s_n \MGL \ar[d]^-{\cov_n} \\
\MGL \ar[r]^-{\tr^n} & f^n \MGL. \\
}
\end{equation}
Here, $u \colon \unit \to \MM\Z$ denotes the unit map.

Recall the following facts about the slices of $\MGL$.

\begin{theorem}\label{thm:SlicesMGL}
\begin{enumerate}
\item The unit map $u \colon \unit \to \MGL$ induces an isomorphism on zeroth slices
\[
s_0(\unit) \ral{\cong} s_0(\MGL) \cong \MM\Z
\]
\cite{Spitzweck10}*{Corollary~3.3}.
\item \emph{If} the Hopkins--Morel--Hoyois Conjecture~\ref{conj:HopkinsMorel} holds over the base scheme $S$, then the slices of $\MGL$ are given by
\begin{align*}
s_n(\MGL) &\cong \Si_{\PP^1}^{n} s_0(\MGL) \ot_{\Z} L_{n} \\
&\cong \Si^{2n,n} \MM\Z \ot_{\Z} L_{n}
\end{align*}
for $n \geq 0$ and trivial for $n < 0$ \cite{Spitzweck10}*{Theorem~4.7}. 
In other words, $S$ satisfies Conjecture~\ref{conj:SliceMGL}. 
\end{enumerate}
\end{theorem}

\begin{theorem}
\cite{Hoyois15}*{Theorem~8.5} Assume that the base scheme $S$ is essentially smooth over a field of characteristic $p > 0$. Then the map %
$\psi_{S} \colon \Si^{2n,n} \MM\Z_{S} \ot_{\Z} L_{n} \to s_n \MGL_{S}$ in $\hMod{\MM\Z_S}$ induces an isomorphism after inverting $p$:
\[
\xymatrix{
\Si^{2n,n} \MM\Z_{S} \ot_{\Z} L_{n} [1/p] \ar[r]^-{\psi_S [1/p]}_-{\cong} & s_n \MGL_{S} [1/p]. \\
}
\]
\end{theorem}

\begin{remark} 
In light of the natural isomorphisms in $\hMod{\MM\Z}$
\begin{align*}
&( \MM\Z ) [1/p] \cong \MM ( \Z [1/p] ) \\
&( s_n E ) [1/p] \cong s_n ( E [1/p] ),  
\end{align*}
there is no ambiguity in the notation. See the discussion after~\cite{Hoyois15}*{Proposition~4.13}.
\end{remark}

We will show a weaker statement than Conjecture~\ref{conj:SliceMGL}, namely:

\begin{theorem}\label{thm:RetractSlice}
Assume that the base scheme $S$ essentially smooth over a field of characteristic $p > 0$. Then the map $\psi_{S}$ in $\hMod{\MM\Z_S}$ is a split monomorphism, i.e., admits a retraction.
\end{theorem}

\begin{proof}
Let $f \colon S \to \Spec(\kk)$ be essentially smooth. Then we have the isomorphism $f^* \MM\Z_{\kk} \ral{\cong} \MM\Z_{S}$ (Theorem~\ref{thm:PullbackHZ}), as well as the natural isomorphism $f^* s_n \cong s_n f^*$ for any $n \in \Z$ (Lemma~\ref{lem:PullbackSlice}). By the commutative diagram in $\hMod{\MM\Z_S}$
\[
\xymatrix @R-0.4pc {
f^* (\Si^{2n,n} \MM\Z_{\kk} \ot_{\Z} L_{n}) \ar[r]^-{f^* \psi_{\kk}} & f^* (s_n \MGL_{\kk}) \ar[d]^{\cong} \\
\Si^{2n,n} f^* \MM\Z_{\kk} \ot_{\Z} L_{n} \ar[d]_{\cong} \ar[u]^{\cong} & s_n f^* \MGL_{\kk} \ar[d]^{\cong} \\
\Si^{2n,n} \MM\Z_{S} \ot_{\Z} L_{n} \ar[r]^-{\psi_{S}} & s_n \MGL_{S}, \\
}
\]
it suffices to prove the claim for $S = \kk$, which we do in the rest of this section.
\end{proof}

\subsection{Construction of the retraction}

As in Notation~\ref{nota:RingMaps}, let $i \colon \Spec(\kk) \inj \Spec(D)$ denote the inclusion of the closed point, and $j \colon \Spec(\KK) \inj \Spec(D)$ its open complement. Consider the truncation map in $\SH(\KK)$
\[
\tr^n \colon \MGL_{\KK} \to f^n \MGL_{\KK}
\]
and the composite in $\SH(D)$
\[
\xymatrix{
\MGL_D \ar[r]^-{\eta} & j_* j^* \MGL_D \cong j_* \MGL_{\KK} \ar[r]^-{j_* \tr^n } & j_* ( f^n \MGL_{\KK} ). \\
}
\]
Applying the functor $i^* \colon \SH(D) \to \SH(\kk)$ yields a map
\[
\xymatrix{
\MGL_{\kk} \cong i^* \MGL_{D} \ar[r] & i^* j_* ( f^n \MGL_{\KK} ). \\
}
\]
Applying the $n$\textsuperscript{th} slice yields a map in $\hMod{\MM\Z_\kk}$
\begin{equation}\label{eq:MapSlices}
\xymatrix{
s_n \MGL_{\kk} \ar[r]^-{} & s_n \left( i^* j_* ( f^n \MGL_{\KK} ) \right). \\
}
\end{equation}

\begin{proposition}\label{pr:SliceOpenClosed}
There is an isomorphism in $\hMod{\MM\Z_{\kk}}$ %
\[
s_n \left( i^* j_* ( f^n \MGL_{\KK} ) \right) \cong \Si^{2n,n} \MM\Z_{\kk} \ot_{\Z} L_{n}.
\]
\end{proposition}

\begin{proof}
Consider the slice filtration of $f^n \MGL_{\KK}$
\begin{equation}\label{eq:SliceFiltTrunc}
0 = f_{n+1} (f^n \MGL_{\KK}) \to f_{n} f^n \MGL_{\KK} \to \ldots \to f_{m+1} (f^n \MGL_{\KK}) \to f_m (f^n \MGL_{\KK}) \to \ldots
\end{equation}
whose $m$\textsuperscript{th} filtration quotient is given by
\[
s_m (f^n \MGL_{\KK}) = \begin{cases}
s_m \MGL_{\KK} &\text{for } m \leq n \\
0 &\text{for } m > n. \\
\end{cases}
\]
Since the Hopkins--Morel--Hoyois isomorphism holds over $\KK$, which is a field of characteristic zero, Theorem~\ref{thm:SlicesMGL} yields
\[
s_m \MGL_{\KK} \cong \Si^{2m,m} \MM\Z_{\KK} \ot_{\Z} L_{m}.
\]
Since the functors $j_*$, $i^*$, and $s_n$ are exact, %
applying $s_n i^* j_*$ to the filtration~\eqref{eq:SliceFiltTrunc} yields a filtration of $s_n i^* j_* (f^n \MGL_{\KK})$ whose $m$\textsuperscript{th} filtration quotient is
\begin{equation}\label{eq:FiltQuotient}
s_n i^* j_* \left( \Si^{2m,m} \MM\Z_{\KK} \ot_{\Z} L_{m} \right)
\end{equation}
for $0 \leq m \leq n$ and zero otherwise. Note that this filtration is still finite. %
By (a slight variant of) Lemma~\ref{lem:EMopenClosed}, there is a splitting
\[
i^* j_* \MM\Z_{\KK} \cong \MM\Z_{\kk} \op \Si^{-1,-1} \MM\Z_{\kk}
\]
in $\SH(\kk)$. This in turn yields
\[
i^* j_* \left( \Si^{2m,m} \MM\Z_{\KK} \ot_{\Z} L_{m} \right) \cong \Si^{2m,m} (\MM\Z_{\kk} \op \Si^{-1,-1} \MM\Z_{\kk}) \ot_{\Z} L_{m},
\]
which has trivial slices above degree $m$, by Theorem~\ref{thm:SlicesHZ} and Lemma~\ref{lem:SliceShift}. Hence, the filtration quotient~\eqref{eq:FiltQuotient} is trivial for $m < n$, and only the $n$\textsuperscript{th} is non-trivial. This shows the isomorphism
\begin{align*}
s_n i^* j_* (f^n \MGL_{\KK}) &\cong s_n i^* j_* \left( \Si^{2n,n} \MM\Z_{\KK} \ot_{\Z} L_{n} \right) \quad \text{(the } n^{\text{th}} \text{ filtration quotient)} \\
&\cong s_n \left( \Si^{2n,n} (\MM\Z_{\kk} \op \Si^{-1,-1} \MM\Z_{\kk}) \ot_{\Z} L_{n} \right) \\
&\cong \Si^{2n,n} s_0 \left( (\MM\Z_{\kk} \op \Si^{-1,-1} \MM\Z_{\kk}) \ot_{\Z} L_{n} \right) \\
&\cong \Si^{2n,n} \left( \MM\Z_{\kk} \ot_{\Z} L_{n} \right)
\end{align*}
as claimed.
\end{proof}

\begin{notation}\label{nota:Retraction}
Let
\[
\te \colon \Si^{2n,n} \MM\Z_{\kk} \ot_{\Z} L_{n} \ral{\cong} s_n \left( i^* j_* ( f^n \MGL_{\KK} ) \right)
\]
denote the isomorphism in $\hMod{\MM\Z_{\kk}}$ from Proposition~\ref{pr:SliceOpenClosed}. Let
\begin{equation}\label{eq:MapSlicesNicer}
\phy \colon s_n \MGL_{\kk} \to \Si^{2n,n} \MM\Z_{\kk} \ot_{\Z} L_{n}
\end{equation}
denote the composite in $\hMod{\MM\Z_{\kk}}$ of the map~\eqref{eq:MapSlices} and the isomorphism $\te^{-1}$, i.e., the composite
\[
\xymatrix @C=\bigcol {
s_n \MGL_{\kk} \ar@/_2.5pc/[rrr]_{\phy} \ar[r]^-{s_n i^* \eta} & s_n i^* j_* \MGL_{{\KK}} \ar[r]^{s_n i^* j_* \tr^n} & s_n i^* j_* ( f^n \MGL_{\KK} ) & \ar[l]_{\te}^{\cong} \Si^{2n,n} \MM\Z_{\kk} \ot_{\Z} L_{n}. \\
}
\]
\end{notation}

\subsection{Proof of Theorem~\ref{thm:RetractSlice}}

We will show that the map $\phy$ is a retraction to $\psi_{\kk}$, i.e., the composite
\[
\xymatrix{
\Si^{2n,n} \MM\Z_{\kk} \ot_{\Z} L_{n} \ar[r]^-{\psi_{\kk}} & s_n \MGL_{\kk} \ar[r]^-{\phy} & \Si^{2n,n} \MM\Z_{\kk} \ot_{\Z} L_{n} \\
}
\]
in $\hMod{\MM\Z_{\kk}}$ is the identity, using Notations~\ref{nota:ApproxSliceDM} and \ref{nota:Retraction}.

By Corollary~\ref{cor:EndoHZ}, 
taking the homotopy group $\pi_{2n,n}$ induces an isomorphism of rings
\[
\xymatrix{
[\Si^{2n,n} \MM\Z_{\kk} \ot_{\Z} L_{n}, \Si^{2n,n} \MM\Z_{\kk} \ot_{\Z} L_{n}] \ar[r]^-{\pi_{2n,n}}_-{\cong} & \Hom_{\Z}(L_{n},L_{n}), \\
}
\]
Hence, it suffices to show that the composite $\phy \circ \psi_{\kk}$ induces the identity on $\pi_{2n,n}$:
\[
\xymatrix @C=\bigcol {
L_{n} \cong \pi_{2n,n} \left( \Si^{2n,n} \MM\Z_{\kk} \ot_{\Z} L_{n} \right) \ar[r]^-{\pi_{2n,n} (\phy \circ \psi_{\kk})}_-{\id ?} & \pi_{2n,n} \left( \Si^{2n,n} \MM\Z_{\kk} \ot_{\Z} L_{n} \right) \cong L_{n}. \\
}
\]
A closer look at the proof of Proposition~\ref{pr:SliceOpenClosed} yields a more precise description of the isomorphism $\te$. Namely, start with the diagram in $\SH(D)$
\[
\xymatrix{
\Si^{2n,n} \MM\Z_{D} \ot_{\Z} L_{n} \ar[r]^-{\eta} & \Si^{2n,n} j_* \MM\Z_{\KK} \ot_{\Z} L_{n} \ar[r]^-{j_* \psi_{\KK}}_-{\cong} & j_* s_n \MGL_{\KK} \ar[r]^{j_* \cov_n} & j_* f^n \MGL_{\KK}.  
}
\]
Applying $s_n i^*$ yields a diagram of isomorphisms in $\hMod{\MM\Z_{\kk}}$
\[
\hspace{-0.05\textwidth}
\resizebox{1.1\textwidth}{!}{
\xymatrix @C=3.3pc {
s_n \left( \Si^{2n,n} \MM\Z_{\kk} \ot_{\Z} L_{n} \right) \ar[r]^-{s_n i^* \eta}_-{\cong} & s_n i^* j_* \left( \Si^{2n,n} \MM\Z_{\KK} \ot_{\Z} L_{n} \right) \ar[r]^-{s_n i^* j_* \psi_{\KK}}_-{\cong} & s_n i^* j_* s_n \MGL_{\KK} \ar[r]^{s_n i^* j_* \cov_n}_-{\cong} & s_n i^* j_* f^n \MGL_{\KK} \\
\Si^{2n,n} \MM\Z_{\kk} \ot_{\Z} L_{n} \ar@{=}[u] \ar@/_1pc/[urrr]_-{\te} & & & \\
}
}
\]
whose composite we called $\te$.

Consider the map $\Si^{2n,n} \unit_{D} \ot_{\Z} L_{n} \to j_* \MGL_{\KK}$ in $\SH(D)$ defined by either composite in the commutative diagram
\[
\xymatrix @-0.4pc {
\Si^{2n,n} \unit_{D} \ot_{\Z} L_{n} \ar[d]_{\eta} \ar[r]^-{\tild{\la}_D} & \MGL_D \ar[d]^{\eta} \\
j_* j^* \left( \Si^{2n,n} \unit_{D} \ot_{\Z} L_{n} \right) \ar[d]_{\cong} \ar[r]^-{j_* j^* \tild{\la}_D} & j_* j^* \MGL_D \ar[d]^{\cong} \\
j_* \left( \Si^{2n,n} \unit_{\KK} \ot_{\Z} L_{n} \right) \ar[r] & j_* \MGL_{\KK}. \\
}
\]
Applying $i^*$ yields a map in $\SH(\kk)$ which we denote
\[
\tild{\ga} \dfn i^*(\eta \circ \tild{\la}_D) \colon \Si^{2n,n} \unit_{\kk} \ot_{\Z} L_{n} \to i^* j_* \MGL_{\KK}.
\]
Let $\ga \colon L_{n} \to \pi_{2n,n} i^* j_* \MGL_{\KK}$ denote the corresponding map of abelian groups.

\begin{lemma}\label{lem:EquivSliceCompat}
The map $\te$ makes the following diagram of abelian groups commute:
\[
\xymatrix @C=\bigcol {
L_{n} \ar[d]_-{\ga} \ar[rr]^-{\text{unit}}_-{\cong} & & \pi_{2n,n} (\Si^{2n,n} \MM\Z_{\kk} \ot_{\Z} L_{n}) \ar[d]^-{\pi_{2n,n} (\te)}_{\cong} \\
\pi_{2n,n} (i^* j_* \MGL_{\KK}) \ar[r]^-{\pi_{2n,n} (i^* j_* \tr^n)} & \pi_{2n,n} (i^* j_* f^n \MGL_{\KK}) \ar[r]^-{\si} & \pi_{2n,n} (s_n i^* j_* f^n \MGL_{\KK}). \\
}
\]
\end{lemma}

\begin{proof}
By adjunction, the statement is equivalent to commutativity of the following diagram in $\SH(\kk)$:
\[
\xymatrix @R-0.4pc {
\Si^{2n,n} \unit_{\kk} \ot_{\Z} L_{n} \ar[dd]_-{\tild{\ga}} \ar[rr]^-{u_{\kk}} & & \Si^{2n,n} \MM\Z_{\kk} \ot_{\Z} L_{n} \ar[d]^-{\te}_{\cong} \\
& & s_n (i^* j_* f^n \MGL_{\KK}) \ar[d]^-{\cov_n} \\
i^* j_* \MGL \ar[r]^-{i^* j_* \tr^n} & i^* j_* f^n \MGL_{\KK} \ar[r]^-{\tr^n} & f^n (i^* j_* f^n \MGL_{\KK}). \\
}
\]
Using the isomorphism
\[
s_n (\Si^{2n,n} \MM\Z_{\kk} \ot_{\Z} L_{n}) = \Si^{2n,n} \MM\Z_{\kk} \ot_{\Z} L_{n}
\]
and the definition of $\te$, this in turn is equivalent to commutativity of the following diagram in $\SH(\kk)$:
\[
\xymatrix @R-0.4pc {
\Si^{2n,n} \unit_{\kk} \ot_{\Z} L_{n} \ar[ddd]_-{\tild{\ga}} \ar[r]^-{u_{\kk}} & \Si^{2n,n} \MM\Z_{\kk} \ot_{\Z} L_{n} \ar[d]^-{i^* \eta} \\
& i^* j_* (\Si^{2n,n} \MM\Z_{\KK} \ot_{\Z} L_{n}) \ar[d]^-{i^* j_* (\psi_{\KK})} \\
& i^* j_* s_n \MGL_{\KK} \ar[d]^-{i^* j_* (\cov_n)} \\
i^* j_* \MGL_{\KK} \ar[r]^-{i^* j_* \tr^n} & i^* j_* f^n \MGL_{\KK}. \\
}
\]
Recalling the definition $\tild{\ga} = i^* (\eta \circ \tild{\la}_D)$, this diagram is $i^*$ applied to the following (outer) diagram in $\SH(D)$:
\[
\xymatrix @R-0.4pc @C=\bigcol {
\Si^{2n,n} \unit_{D} \ot_{\Z} L_{n} \ar[d]_-{\eta} \ar[r]^-{u_{D}} & \Si^{2n,n} \MM\Z_{D} \ot_{\Z} L_{n} \ar[d]^-{\eta} \\
j_*  \left( \Si^{2n,n} \unit_{\KK} \ot_{\Z} L_{n} \right) \ar[dd]_-{j_* j^* \tild{\la}_D = j_* \tild{\la}_{\KK}} \ar[r]^-{j_* j^* u_D = j_* u_{\KK}} & j_* (\Si^{2n,n} \MM\Z_{\KK} \ot_{\Z} L_{n}) \ar[d]^-{j_* (\psi_{\KK})}_{\cong} \\
& j_* s_n \MGL_{\KK} \ar[d]^-{j_* (\cov_n)} \\
j_* \MGL_{\KK} \ar[r]^-{j_* \tr^n} & j_* f^n \MGL_{\KK}. \\
}
\]
Note that we used the equation $j^* \tild{\la}_D = \tild{\la}_{\KK}$ from Lemma~\ref{lem:LazardMap}. The top part of the diagram commutes by naturality of $\eta$. The bottom part commutes because it is $j_*$ applied to the commutative diagram~\eqref{eq:CompatLazardSH}.%
\end{proof}

Now consider the commutative diagram of abelian groups
\[
\hspace{-0.05\textwidth}
\resizebox{1.1\textwidth}{!}{
\xymatrix{
\pi_{2n,n} (\Si^{2n,n} \MM\Z_{\kk} \ot_{\Z} L_{n}) \ar@/_1.5pc/[ddr]_-{\pi_{2n,n}(\psi_{\kk})} & L_{n} \ar[l]_-{\text{unit}}^-{\cong} \ar[d]_-{\la_{\kk}} \ar@/^1pc/[dr]^{\ga} & & \\
& \pi_{2n,n} \MGL_{\kk} \ar[r]^-{\pi_{2n,n}(i^* \eta)} \ar[d]_{\si} & \pi_{2n,n} i^* j_* \MGL_{\KK} \ar[r]^-{\pi_{2n,n}(i^* j_* \tr^n)} & \pi_{2n,n} i^* j_* f^n \MGL_{\KK} \ar[d]^{\si} \\
& \pi_{2n,n} (s_n \MGL_{\kk}) \ar[drr]_{\pi_{2n,n} (\phy)} \ar[rr] & & \pi_{2n,n} (s_n i^* j_* f^n \MGL_{\KK}) \\
& & & \pi_{2n,n} (\Si^{2n,n} \MM\Z_{\kk} \ot_{\Z} L_{n}) \ar[u]_{\pi_{2n,n} (\te)}^{\cong} \\
}
}
\]
where the middle rectangle comes from naturality of $\si$. The left triangle commutes by definition of $\psi_{\kk}$, while the bottom right triangle commutes by definition of $\phy$. %
The top triangle commutes by definition of $\ga$. Combined with Lemma~\ref{lem:EquivSliceCompat}, this diagram yields the equation $\pi_{2n,n}(\phy \circ \psi_{\kk}) = \id$, 
completing the proof.

\section{Reduction to smoothness of Eilenberg--MacLane spectra}\label{sec:ReduceLisse}

\begin{definition}\label{def:Smooth}
An object of $\SH(S)$ is \Def{lisse} 
if it lies in the full localizing triangulated subcategory of $\SH(S)$ generated by the strongly dualizable motivic spectra.
\end{definition}

\begin{remark}
In the (classical) stable homotopy category, every spectrum in the full localizing triangulated subcategory generated by the sphere spectrum $S^0$, which is strongly dualizable, with dual $S^0$. Therefore, every object of the stable homotopy category is lisse in the sense of Definition~\ref{def:Smooth}. This is the case in the homotopy category of several familiar stable homotopy theories, as discussed in \cite{Strickland04}*{\S 1}.

For motivic spectra, R\"ondigs and {\O}stv{\ae}r showed that every object of $\SH(\KK)$ is lisse if $\KK$ is a field of characteristic zero, using resolution of singularities \cite{RondigsO08}*{Theorem~52, Remark~59}. See also the discussion at \cite{Hoyois15MO}.
\end{remark}

\begin{lemma}\label{lem:CompactnessCoprod}
Let $F \colon \cat{C} \rla \cat{D} \colon G$ be an adjoint pair of triangulated functors between triangulated categories with arbitrary coproducts. Assume that $\cat{C}$ is compactly generated. Then the following are equivalent.
\begin{enumerate}
\item The right adjoint $G$ preserves coproducts.
\item The left adjoint $F$ preserves compactness.
\item For some set of compact generators $\{ C_{\al} \}$ of $\cat{C}$, every object $F(C_{\al})$ is compact in $\cat{D}$. 
\end{enumerate}
\end{lemma}

\begin{proof}
See \cite{Neeman96}*{Theorem~5.1}.
\end{proof}

\begin{lemma}\label{lem:PreserveCoprod}
For any map of base schemes $f \colon S \to T$, the functor $f^* \colon \SH(T) \to \SH(S)$ preserves compact objects. 
\end{lemma}

As a consequence, the right adjoint $f_* \colon \SH(S) \to \SH(T)$ preserves coproducts.

\begin{proof}
Since $S$ is Noetherian and of finite Krull dimension, the category $\Sm_S$ is essentially small; see \cite{DundasRO03}*{\S 2}. Via the Yoneda embedding, view an object $X$ of $\Sm_S$ as a motivic space over $S$. Its suspension spectrum $\Si^{\infty} X_+$ is a compact object of $\SH(S)$. As $X$ ranges over isomorphism classes of objects in $\Sm_S$, the objects $\Si^{0,q} \Si^{\infty} X_+$ with integers $q \leq 0$ 
form a generating set of $\SH(S)$; see \cite{DuggerI05}*{Proposition~9.2} or \cite{Morel04int}*{Example~5.8, Proposition~5.1.14}. 
Now, for $Y$ in $\Sm_T$, the pullback
\[
f^* \left( \Si^{0,q} \Si^{\infty} Y_+ \right) = \Si^{0,q} \Si^{\infty} f^*Y_+
\]
is compact in $\SH(S)$. 
By Lemma~\ref{lem:CompactnessCoprod}, $f^*$ preserves compactness.
\end{proof}

\begin{notation}
Let $f \colon S \to T$ be a map of base schemes. For $Y$ in $\SH(T)$ and $X$ in $\SH(S)$, consider the \Def{exchange transformation} in $\SH(T)$
\[
\exch{f} \colon f_* X \sm_T Y \to f_* \left( X \sm_S f^* Y \right)
\]
which is the composite
\[
\xymatrix{
f_* X \sm_T Y \ar[r]^-{\eta} & f_* f^* \left( f_* X \sm_T Y \right) \ar@{-}[r]^-{\cong} & f_* \left( f^* f_* X \sm_S f^* Y \right) \ar[r]^-{f_*(\ep \sm 1)} & f_* \left( X \sm_S f^* Y \right), \\
}
\]
cf.\ \cite{CisinskiD19}*{1.1.31}.
\end{notation}

\begin{lemma}\label{lem:ProjFormula}
If $Y$ is a lisse object of $\SH(T)$, then $\exch{f}$ is an isomorphism. In other words, $Y$ and $f_*$ satisfy the \emph{projection formula} $f_* X \sm_T Y \cong f_* \left( X \sm_S f^* Y \right)$.
\end{lemma}

\begin{proof}
First, note that both sides of the %
exchange transformation 
$\exch{f}$ commute with (homotopy) colimits in the variable $Y$. Indeed, the functors $f_* X \sm_T -$, $f^*$, and $X \sm_S -$ are exact and left adjoints, and thus preserve homotopy colimits. The functor $f_*$ is exact and preserves coproducts, by Lemma~\ref{lem:PreserveCoprod}. Since $Y$ is obtained as a sequential homotopy colimit
\[
Y = \hocolim \left( 0 \to Y_0 \to Y_1 \to Y_2 \to \ldots \right)
\]
where the cofiber of each $Y_i \to Y_{i+1}$ is a coproduct of suspensions of strongly dualizable objects \cite{HoveyPS97}*{Proposition 2.3.17}, it suffices to prove the claim for $Y$ strongly dualizable. This last step is done in \cite{FauskHM03}*{Proposition~3.2}, 
which applies here since $\SH(S)$ is symmetric monoidal and $f^* \colon \SH(T) \to \SH(S)$ is strong monoidal.
\end{proof}

\begin{remark}\label{rem:PushforwardTwist}
Taking $Y = S^{0,1}$ in Lemma~\ref{lem:ProjFormula} shows that a direct image functor $f_*$ commutes with twists (or equivalently, with $\G_m$-suspension). This yields a natural isomorphism $f_*(\Si^{p,q} X) \cong \Si^{p,q} f_* X$ in $\SH(T)$. 
\end{remark}

A straightforward diagram chase yields the following.

\begin{lemma}\label{lem:ExchangeSmash}
Let $f \colon S \to T$ be a map of base schemes. For any $X, Y$ in $\SH(T)$, the exchange transformation $\exch{f}$ makes the following diagram in $\SH(T)$ commute:
\[
\xymatrix @C+0.3pc {
X \sm Y \ar[r]^-{\eta \sm 1} \ar[drr]_-{\eta} & f_* f^* X \sm Y \ar[r]^-{\exch{f}} & f_* (f^* X \sm f^* Y) \ar[d]^{\cong} \\ 
& & f_* f^* (X \sm Y). \\
}
\]
\end{lemma}

\begin{lemma}\label{lem:ShriekCoprod}
Let $i \colon Z \to T$ be a closed immersion of base schemes. Then the functor $i^! \colon \SH(T) \to \SH(Z)$ preserves coproducts.
\end{lemma}

\begin{proof}
The functor $i^!$ sits in the exact triangle~\eqref{eq:ExactTriK}, where the functors $i^*$, $j_*$, and $j^*$ preserve coproducts; see \cite{Ayoub07I}*{Lemme~2.1.157}.
\end{proof}

\begin{notation}\label{nota:ExchangeModule}
Let $i \colon Z \to T$ be a closed immersion of base schemes. For any $X, Y$ in $\SH(T)$, denote by
\[
\exmodr \colon i^! X \sm i^* Y \to i^!(X \sm Y)
\]
the natural map in $\SH(Z)$ described in \cite{Ayoub07I}*{\S 2.3.2}. %
Explicitly, $\exmodr$ is the composite 
\begin{equation}\label{eq:Exmod}
\xymatrix @C+0.5pc {
i^! X \sm i^* Y \ar[r]^-{\eta} & i^! i_* (i^! X \sm i^* Y) & i^! (i_* i^! X \sm Y) \ar[l]_-{i^! \exch{i}}^-{\cong} \ar[r]^-{i^! (\ep \sm 1)} & i^!(X \sm Y), \\
}
\end{equation}
where $\ep \colon i_* i^! \to 1$ denotes the counit of the adjunction $i_* \dashv i^!$. The exchange transformation $\exch{i}$ is an isomorphism for $i$ a closed immersion \cite{Ayoub07I}*{Lemme~2.3.10}. 
The map $\exmodl \colon i^* Y \sm i^! X \to i^!(Y \sm X)$ is defined similarly.
\end{notation}

Equation~\eqref{eq:Exmod} can be reformulated as commutativity of the following diagram in $\SH(T)$:
\[
\xymatrix{
i_* (i^! X \sm i^* Y) \ar[r]^-{i_* \exmodr} & i_* i^! (X \sm Y) \ar[r]^-{\ep} & X \sm Y \\ 
i_* i^! X \sm Y. \ar[u]^{\exch{i}}_{\cong} \ar[urr]_-{\ep \sm 1} & & \\
}
\]
Applying $i^*$ yields the commutative diagram in $\SH(Z)$:
\[
\xymatrix{
i^! X \sm i^* Y \ar[r]^-{\exmodr} \ar[drr]_-{i^* \ep \sm 1} & i^! (X \sm Y) \ar[r]^-{i^* \ep} & i^* (X \sm Y) \\ 
& & i^* X \sm i^* Y. \ar[u]_{\cong} \\
}
\]
Combined with Lemma~\ref{lem:ExchangeSmash}, this yields the following.

\begin{lemma}\label{lem:ExchangeTriangle}
Let $i \colon Z \to T$ be a closed immersion of base schemes. 
For any $X, Y$ in $\SH(T)$,  the following diagram in $\SH(Z)$ commutes:
\[
\xymatrix{
i^! X \sm i^* Y \ar[d]^{\exmodr} \ar[r]^-{i^* \ep \sm 1} & i^*X \sm i^*Y \ar[d]^{\cong} \ar[r]^-{i^* \eta \sm 1} & i^* j_* j^* X \sm i^* Y \ar[d]^{i^* \exch{j}} \\
i^! (X \sm Y) \ar[r]^-{i^* \ep} & i^* (X \sm Y) \ar[r]^-{i^* \eta} & i^* j_* j^* (X \sm Y). \\
}
\]
The top row is the exact triangle~\eqref{eq:ExactTriK} for $X$ smashed with $i^* Y$,  whereas the bottom row is the triangle for $X \sm Y$.
\end{lemma}

\begin{proposition}\label{pr:ReduceLisse}
Assume that the map in $\SH(\kk)$
\[
i^* \exch{j} \colon (i^* j_* \MM\F_p^{\KK}) \sm \MM\F_p^{\kk} \to i^* j_* (\MM\F_p^{\KK} \sm \MM\F_p^{\KK})
\]
is an isomorphism. 
Then the map $\Psi^{\kk}$ in $\hMod{\MM\F_p^{\kk}}$ is an isomorphism. In other words, the splitting of the dual Steenrod algebra given in Theorem~\ref{thm:DualSteenrod} holds over the base scheme $S = \Spec(\kk)$.
\end{proposition}

The assumption holds for instance if $\hatt{M}\F_{p}^{D}$ is lisse in $\SH(D)$, by Lemma~\ref{lem:ProjFormula}.

\begin{proof}
We will show that the underlying map of spectra $\Psi^{\kk}$ in $\SH(\kk)$ is an isomorphism. Pick a cofiber of $\Psi^{D}$, i.e., an exact triangle in $\SH(D)$
\[
\xymatrix{
\bigoplus_{\seq} \Si^{p_{\seq},q_{\seq}} \MM\F_{p}^D \ar[r]^-{\Psi^D} & \MM\F_{p}^D \sm \MM\F_{p}^D \ar[r]^-{q} & \Cof(\Psi^D). \\
}
\]
Applying $i^*$ yields an exact triangle in $\SH(\kk)$, which we put as middle column of the diagram:
\begin{equation}\label{eq:Cofibers}
\xymatrix @C=\bigcol {
i^! (\bigoplus_{\seq} \Si^{p_{\seq},q_{\seq}} \MM\F_{p}^D) \ar@{^{(}->}[d]^{i^! \Psi^D} \ar[r]^-{i^* \ep = 0} & \bigoplus_{\seq} \Si^{p_{\seq},q_{\seq}} \MM\F_{p}^{\kk} \ar@{^{(}->}[d]^{\Psi^{\kk}} \ar@{^{(}->}[r]^-{i^* \eta} & \bigoplus_{\seq} \Si^{p_{\seq},q_{\seq}} i^* j_* \MM\F_{p}^{\KK} \ar[d]^{i^* j_* \Psi^{\KK}}_{\simeq} \\
i^! (\MM\F_p^{D} \sm \MM\F_p^{D}) \ar@{->>}[d]^{i^! q} \ar[r]^-{i^* \ep=0} & \MM\F_p^{\kk} \sm \MM\F_p^{\kk} \ar[d]^{i^* q} \ar[r]^-{i^* \eta} & i^* j_* (\MM\F_p^{\KK} \sm \MM\F_p^{\KK}) \ar[d]^{i^* j_* j^* q} \\
i^! \Cof(\Psi^D) \ar[r]^-{i^* \ep}_-{\therefore 0} & i^* \Cof(\Psi^D) \ar[r]^-{i^* \eta} & i^* j_* j^* \Cof(\Psi^D) = 0. \\
}
\end{equation}
Here, the rows are the natural localization triangle~\eqref{eq:ExactTriK}. 
Note that each column is exact, since the functors $i^!$, $i^*$, $j_*$, and $j^*$ are exact. In the top row, the map $i^* \ep$ is zero, since the map $i^* \ep \colon i^! \MM\F_p^D \to \MM\F_p^{\kk}$ is zero, and $i^!$ preserves coproducts (Lemma~\ref{lem:ShriekCoprod}).

Now we show that in the middle row, the map $i^* \ep$ is zero. Lemma~\ref{lem:ExchangeTriangle} provides a map of exact triangles
\[
\xymatrix @C=\bigcol {
(i^! \MM\F_p^{D}) \sm \MM\F_p^{\kk} \ar[d]^{\exmodr}_{\therefore \cong} \ar[r]^-{i^* \ep \sm 1}_-{=0 \sm 1} & \MM\F_p^{\kk} \sm \MM\F_p^{\kk} \ar@{=}[d] \ar[r]^-{i^* \eta \sm 1} & (i^* j_* \MM\F_p^{\KK}) \sm \MM\F_p^{\kk} \ar[d]^{i^* \exch{j}}_{\cong} \\
i^! (\MM\F_p^{D} \sm \MM\F_p^{D}) \ar[r]^-{i^* \ep}_-{\therefore 0} & \MM\F_p^{\kk} \sm \MM\F_p^{\kk} \ar[r]^-{i^* \eta} & i^* j_* (\MM\F_p^{\KK} \sm \MM\F_p^{\KK}) \\
}
\]
which is an isomorphism of triangles, by the assumption that $i^* \exch{j}$ is an isomorphism. Moreover, since $\exmodr \colon i^! \MM\F_{p}^{D} \sm \MM\F_{p}^{\kk} \to i^!(\MM\F_{p}^D \sm \MM\F_{p}^D)$ is an isomorphism, Corollary~\ref{cor:ShriekSteenrod} identifies the map $i^! \Psi^{D}$ with $\Si^{-2,-1} \Psi^{\kk}$, which is a split monomorphism (Theorem~\ref{thm:RetractSteenrod}).

Finally, the map $i^* \ep$ in the bottom row of Diagram~\eqref{eq:Cofibers} is zero, since $i^! q$ is a (split) epimorphism. This implies $i^* \Cof(\Psi^D) = 0$, so that $\Psi^{\kk}$ is an isomorphism.
\end{proof}

Combining Proposition~\ref{pr:ReduceLisse} and Corollary~\ref{cor:ReduceFp} yields a proof of Proposition~\ref{pr:ReduceLisseIntro}.

\subsection{Reduction to Eilenberg--MacLane spaces}

One approach using Proposition~\ref{pr:ReduceLisse} %
would be to find sufficient conditions on motivic Eilenberg--MacLane \emph{spaces} to guarantee that $\MM\F_p$ is lisse. 
Let us collect a few observations about lisse objects in $\SH(S)$.

\begin{lemma}
\begin{enumerate}
\item \label{item:SmoothProj} If $X$ is a smooth and proper $S$-scheme, viewed as a motivic space over $S$, then its $\PP^1$-suspension spectrum $\Si^{\infty} X_+$ is dualizable in $\SH(S)$, in particular lisse.
\item \label{item:Hocolim} The property of being lisse in $\SH(S)$ is closed under homotopy colimits.
\item \label{item:Suspension} The property of a motivic space $X$ over $S$ that ``its $\PP^1$-suspension spectrum $\Si^{\infty} X_+$ is lisse in $\SH(S)$'' is closed under homotopy colimits. 
\end{enumerate}
\end{lemma}

\begin{proof}
For (\ref*{item:SmoothProj}), see \cite{Hoyois17}*{Corollary~6.13} or \cite{Hu05}*{Appendix~A}. 
For (\ref*{item:Hocolim}), a localizing subcategory is closed under homotopy colimits. 
For (\ref*{item:Suspension}), both the disjoint basepoint functor $(-)_+ \colon \Spc(S) \to \Spc(S)_*$ and suspension spectrum functor $(-)_+ \colon \Spc(S)_* \to \Spt(S)$ preserve homotopy colimits.
\end{proof}

We have the equivalence in $\SH(S)$
\[
\MM\Z = \hocolim\limits_n \Si^{-2n,-n} \Si^{\infty} K(\Z(n),2n),
\]
where $K(\Z(q),p)$ denotes the motivic Eilenberg--MacLane space, also denoted $K(\Z,p,q)$ in the literature. 
Hence, it would suffice to show that $\Si^{\infty} K(\Z(n),2n)$ is lisse for $n$ large enough to conclude that $\MM\Z$ is lisse. 
A similar description holds for $\hatt{M}\Z$.

\section{Purity for motivic cohomology}\label{sec:Purity}

Let $i \colon Z \inj X$ be a regular closed immersion between quasi-compact quasi-separated schemes. Let $N_i$ denote the normal bundle 
of $i \colon Z \inj X$ and $\Thom(N_i)$ its Thom spectrum, as an object of $\SH(Z)$. 
Let $\Thom(-N_i) \cong \Thom(N_i)^{\dual}$ denote the Thom spectrum of the virtual bundle $-N_i$, which is the dual of $\Thom(N_i)$ in $\SH(Z)$. We will use a construction due to D\'eglise, Jin, and Khan of a certain map
\[
\pur \colon i^*A \sm \Thom(N_i)^{\dual} \to i^! A
\]
which is natural in the object $A$ of $\SH(X)$, i.e., a natural transformation of functors $\SH(X) \to \SH(Z)$, called the \emph{purity transformation} \cite{DegliseJK21}.

\begin{lemma}\label{lem:PurityBimod}
The purity transformation is a morphism of $i^*$-bimodules, i.e., makes the following diagrams in $\SH(Z)$ commute:
\[
\xymatrix{
(i^* A \sm \Thom(N_i)^{\dual}) \sm i^* B \ar[d]_{\cong} \ar[r]^-{\pur_A \sm 1} & i^! A \sm i^* B \ar[d]^{\exmodr} \\
i^* (A \sm B) \sm \Thom(N_i)^{\dual} \ar[r]^-{\pur_{A \sm B}} & i^! (A \sm B) \\
} \hspace{1.5cm}
\xymatrix{
i^* B  \sm (i^* A \sm \Thom(N_i)^{\dual}) \ar[d]_{\cong} \ar[r]^-{1 \sm \pur_A} & i^* B \sm i^! A  \ar[d]^{\exmodl} \\
i^* (B \sm A) \sm \Thom(N_i)^{\dual} \ar[r]^-{\pur_{B \sm A}} & i^! (B \sm A). \\
}
\]
\end{lemma}

\begin{proof}
The following argument was kindly provided to us by Fr\'ed\'eric D\'eglise. There is a certain ``fundamental class'' 
$c \colon \Thom(N_i)^{\dual} \to i^!(\unit)$ 
and the purity transformation agrees with the composite
\[
\xymatrix{
i^*A \sm \Thom(N_i)^{\dual} \ar[r]^-{1 \sm c} & i^*A \sm i^!(\unit) \ar[r]^-{\exmodl} & i^! (A \sm \unit) = i^! A.
}
\]
That $\pur$ commutes with the right and left $i^*$-module structures follows from the associativity equations for the $i^*$-bimodule $i^!$ \cite{DegliseJK21}*{2.1.10}.
\end{proof}

\begin{corollary}
Let $R$ be a motivic $E_{\infty}$ ring spectrum which is Cartesian. Then the purity transformation $\pur \colon i^*(-) \sm \Thom(N_i)^{\dual} \to i^!$ lifts canonically to a natural transformation of functors $\hMod{R_X} \to \hMod{R_Z}$.
\end{corollary}

For the remainder of the section, specialize to our example of interest, the closed immersion $i \colon \Spec(\kk) \to \Spec(D)$ from Notation~\ref{nota:RingMaps}. In that case, the normal bundle $N_i$ is a trivial line bundle over $\Spec(\kk)$, whose Thom space is the cofiber of the map $\A^1_{\kk} \setminus \{0\} \to \A^1_{\kk}$. This yields an equivalence $\Thom(N_i) \simeq \PP^1_{\kk} \simeq S^{2,1}$ of motivic spaces and also their suspension spectra, %
so that the purity transformation has the form $\pur \colon \Si^{-2,-1} i^* \to i^!$.

\begin{proposition}\label{pr:ShriekHZ}
Evaluated at the object $\hatt{M}\Z_D$ of $\SH(D)$, the purity transformation in $\SH(\kk)$
\[
\pur_{\hatt{M}\Z} \colon \Si^{-2,-1} i^* \hatt{M}\Z_{D} \ral{\cong} i^! \hatt{M}\Z_{D} 
\]
is an isomorphism.
\end{proposition}

\begin{corollary}\label{cor:ShriekHZ}
The purity transformation provides an isomorphism in $\hMod{\MM\Z_{\kk}}$
\[
i^! \MM\Z_{D} \cong \Si^{-2,-1} \MM\Z_{\kk}
\]
and an isomorphism in $\hMod{\MM\F_{\ell}^{\kk}}$ 
\[
i^! \MM\F_{\ell}^D \cong \Si^{-2,-1} \MM\F_{\ell}^{\kk}.
\]
\end{corollary}

\begin{proof}[of Proposition~\ref{pr:ShriekHZ}]
By %
\cite{Spitzweck18}*{Theorem~7.4}, 
there is an isomorphism $i^! \hatt{M}\Z_{D} \cong \Si^{-2,-1} \MM\Z_{\kk}$ in $\SH(\kk)$. 
By Corollary~\ref{cor:EndoHZ}, the map $\pur_{\hatt{M}\Z}$ is an isomorphism if and only if it induces a surjection on the homotopy group $\pi_{-2,-1}$. Since %
$\Si^{-2,-1} \MM\Z_{\kk}$ is %
its own $(-1)$-slice, the latter condition is equivalent to $s_{-1} \pur_{\hatt{M}\Z}$ inducing a surjection on $\pi_{-2,-1}$, which is what we will show.

Consider the composite in $\SH(D)$
\[
\xymatrix{
s_0 \unit_D \ar[r]^-{s_0 u} & s_0 \MM\Z_D & f_0 \MM\Z_D \ar[l]_-{\tr^0}^{\cong} \ar[r]^-{\cov_0} & \MM\Z_D \\ 
}
\]
using the fact that $\MM\Z_{D}$ is $0$-truncated %
\cite{Spitzweck18}*{Remark~10.2}. 
Consider the diagram in $\SH(\kk)$
\[
\xymatrix @C = \bigcol {
\Si^{-2,-1} i^* \MM\Z_{D} \ar[r]^-{\pur_{\MM\Z}} & i^! \MM\Z_D \\
\Si^{-2,-1} i^* s_0 \unit_{D} \ar[u]^{\Si^{-2,-1} i^* (\cov_0 \circ s_0 u)} \ar[d]_{\Si^{-2,-1} i^* (s_0 u)} \ar[r]^-{\pur_{s_0 \unit}} & i^! s_0 \unit_{D} \ar[u]_{i^! (\cov_0 \circ s_0 u)} \ar[d]^{i^! (s_0 u)} \\
\Si^{-2,-1} i^* s_0 \KGL_{D} \ar[r]^-{\pur_{s_0 \KGL}} & i^! s_0 \KGL_D \\
\Si^{-2,-1} i^* f_0 \KGL_{D} \ar[u]^{\Si^{-2,-1} i^* \tr^0} \ar[d]_{\Si^{-2,-1} i^* \cov_0} \ar[r]^-{\pur_{f_0 \KGL}} & i^! f_0 \KGL_{D} \ar[u]_{i^! \tr^0} \ar[d]^{i^! \cov_0} \\
\Si^{-2,-1} i^* \KGL_{D} \ar[r]^-{\pur_{\KGL}}_-{\cong} & i^! \KGL_{D}. \\
}
\]
The purity transformation $\pur_{\KGL}$ is an isomorphism since $\KGL$ satisfies absolute purity %
\cite{CisinskiD19}*{Theorem~13.6.3}. 
Applying the $(-1)$-slice $s_{-1}$ yields the diagram in $\SH(\kk)$
\begin{equation}\label{eq:PurityKGL}
\xymatrix @C=\bigcol {
\Si^{-2,-1} \MM\Z_{\kk} \ar[r]^-{\cong} & \Si^{-2,-1} s_0 i^* \MM\Z_{D} \ar[r]^-{s_{-1} \pur_{\MM\Z}} & s_{-1} i^! \MM\Z_D \\
& \Si^{-2,-1} s_0 i^* s_0 \unit_{D} \ar[u]^{\Si^{-2,-1} s_0 i^* (\cov_0 \circ s_0 u)} \ar[d]_{\Si^{-2,-1} s_0 i^* (s_0 u)} \ar[r]^-{s_{-1} \pur_{s_0 \unit}} & s_{-1} i^! s_0 \unit_{D} \ar[u]_{s_{-1} i^! (\cov_0 \circ s_0 u)}^{\cong} \ar[d]^{s_{-1} i^! (s_0 u)}_{\cong} \\
& \Si^{-2,-1} s_0 i^* s_0 \KGL_{D} \ar[r]^-{s_{-1} \pur_{s_0 \KGL}} & s_{-1} i^! s_0 \KGL_D \\
& \Si^{-2,-1} s_0 i^* f_0 \KGL_{D} \ar[u]^{\Si^{-2,-1} s_0 i^* \tr^0} \ar[d]_{\Si^{-2,-1} s_0 i^* \cov_0} \ar[r]^-{s_{-1} \pur_{f_0 \KGL}} & s_{-1} i^! f_0 \KGL_{D} \ar[u]_{s_{-1} i^! \tr^0}^{\cong} \ar[d]^{s_{-1} i^! \cov_0}_{\cong} \\
\Si^{-2,-1} \MM\Z_{\kk} \ar[r]^-{\cong} & \Si^{-2,-1} s_0 i^* \KGL_{D} \ar[r]^-{s_{-1} \pur_{\KGL}}_-{\cong} & s_{-1} i^! \KGL_{D}. \\
}
\end{equation}
We used the isomorphisms from Lemma~\ref{lem:ShriekCoverKGL} \eqref{item:IsoSlicesShriek} in the right column. 
The map $\Si^{-2,-1} s_0 i^* \cov_0$ in the left column is surjective on $\pi_{-2,-1}$ since $s_0 i^* \cov_0$ is a map of ring spectra and thus preserves the unit. Likewise, $s_0 i^* \tr^0$ and $s_0 i^* (s_0 u)$ preserve the unit. 
The bottom square of Diagram~\eqref{eq:PurityKGL} shows that $s_{-1} i^! \cov_0$ is surjective on $\pi_{-2,-1}$, which implies that it is an isomorphism, by Lemma~\ref{lem:ShriekCoverKGL} \eqref{item:SliceCoverKGL}. %
Hence, the map $s_{-1} \pur_{\MM\Z}$ is also surjective on $\pi_{-2,-1}$.
\end{proof}

\begin{lemma}\label{lem:ShriekCoverKGL}
\begin{enumerate}
\item \label{item:SliceCoverKGL} The object $i^! f_0 \KGL_D$ in $\SH(\kk)$ has $(-1)$-slice
\[
s_{-1} (i^! f_0 \KGL_D) \cong \Si^{-2,-1} \MM\Z_{\kk}.
\]
\item \label{item:IsoSlicesShriek} The three maps in $\SH(\kk)$
\[
\xymatrix @C=\bigcol {
i^! f_0 \KGL_D \ar[r]^-{i^! \tr^0} & i^! s_0 \KGL_D & i^! s_0 \unit_D \ar[l]_-{i^! (s_0 u)} \ar[r]^-{i^! (\cov_0 \circ s_0 u)} & i^! \MM\Z_D \\
}
\]
induce isomorphisms on $(-1)$-slices.
\end{enumerate}
\end{lemma}

\begin{proof}
(1) Consider the localization triangle~\eqref{eq:ExactTriK} for $f_0 \KGL_{D}$. Applying the exact functor $s_{-1}$ yields the exact triangle in $\SH(\kk)$
\[
\xymatrix{
s_{-1} i^! f_0 \KGL_D \ar[r]^-{} & s_{-1} i^* f_0 \KGL_D \ar[r]^-{} & s_{-1} i^* j_* j^* f_0 \KGL_D. \\
}
\]
Since $i^* f_0 E$ is $0$-effective for any object $E$ of $\SH(D)$, we obtain in particular $s_{-1} i^* f_0 \KGL_D = 0$, and thus an isomorphism
\[
s_{-1} i^* j_* j^* f_0 \KGL_D \cong \Si^{1,0} s_{-1} i^! f_0 \KGL_D.
\]
Let us compute the left-hand side. Since $j$ is smooth, we have
\[
i^* j_* j^* f_0 \KGL_D \cong i^* j_* f_0 j^* \KGL_D \cong i^* j_* f_0 \KGL_{\KK}.
\]
The last isomorphism uses the fact that $KGL$ is a Cartesian ring spectrum \cite{CisinskiD19}*{\S 13.1}. %
Starting from the slice filtration of $f_0 \KGL_{\KK}$, applying the exact functor $i^* j_*$ yields a filtration on $i^* j_* f_0 \KGL_{\KK}$ with $m$\textsuperscript{th} filtration quotient 
\[
i^* j_* s_m (f_0 \KGL_{\KK}) = \begin{cases}
i^* j_* (\Si^{2m,m} \MM\Z_{\KK}) &\text{if } m \geq 0, \\
0 &\text{if } m < 0.
\end{cases}
\]
Note that $\Si^{2m,m} i^* j_* \MM\Z_{\KK}$ has slices concentrated in degrees $m$ and $m-1$. Indeed, the localization triangle~\eqref{eq:ExactTriK} for $\MM\Z_{D}$ together with the isomorphism $i^! \MM\Z_{D} \cong \Si^{-2,-1} \MM\Z_{\kk}$ in $\SH(\kk)$ yield
\[
s_m (i^* j_* \MM\Z_{\KK}) = \begin{cases}
\MM\Z_{\kk} &\text{if } m = 0, \\
\Si^{-1,-1} \MM\Z_{\kk} &\text{if } m = -1, \\
0 &\text{otherwise}.
\end{cases}
\]
Therefore, $s_{-1} i^* j_* j^* f_0 \KGL_D$ has a filtration with a single non-zero filtration quotient (in degree $0$)
\[
s_{-1} i^* j_* \MM\Z_{\KK} = \Si^{-1,-1} \MM\Z_{\kk},
\]
which provides the claimed isomorphism, using Lemma~\ref{lem:CompleteFiltration}.

(2) The argument in part (1) also applies to $s_0 \KGL_D$. Naturality of the localization triangle~\eqref{eq:ExactTriK} and of the slice filtration implies that $i^! \tr^0 \colon i^! f_0 \KGL_D \to i^! s_0 \KGL_D$ induces an isomorphism on $(-1)$-slices.

The same argument applies to $s_0 \unit_D$. %
Since the isomorphisms $s_0 \unit_{\KK} \cong \MM\Z_{\KK}$ and $s_0 \KGL_{\KK} \cong \MM\Z_{\KK}$ are induced by the respective units $u \colon \unit_{\KK} \to \MM\Z_{\KK}$ and $u \colon \unit_{\KK} \to \KGL_{\KK}$ (Theorem~\ref{thm:SlicesHZ}), the claimed isomorphisms of $(-1)$-slices follow.
\end{proof}

For the next lemma, denote by $\SH(S)_{\geq n}$ the subcategory of $\SH(S)$ generated under homotopy colimits and extensions by objects $\Si^{p,q} \Si^{\infty} X_+$ for $X \in \Sm_S$ and $p-q \geq n$. Denote the full subcategory of $\SH(S)$
\[
\SH(S)_{h\geq n} = \{ E \in \SH(S) \mid \ul{\pi}_{p,q} E = 0 \text{ for } p-q < n \},
\]
where the $\ul{\pi}_{p,q} E$ denote the homotopy sheaves of $E$. The inclusion $\SH(D)_{h \geq n} \subseteq \SH(D)_{\geq n}$ always holds, and is an equality over a base field %
\cite{Spitzweck20}*{Lemma~4.3} 
\cite{Hoyois15}*{Theorem~2.3}. %
 
\begin{lemma}\label{lem:CompleteFiltration}
In $\SH(\kk)$, we have $\holim\limits_{n \to +\infty} s_{-1} i^* j_* (f_n \KGL_{\KK}) = 0$.
\end{lemma}

\begin{proof}
Since $\KGL_{\KK}$ is slice complete %
\cite{RondigsSO18}*{Lemma~3.11}, 
we have
\begin{align*}
j_* (f_n \KGL_{\KK}) &= j_* (\holim_{m \to +\infty} f^m f_n \KGL_{\KK}) \\
&= \holim_{m \to +\infty} j_* (f^m f_n \KGL_{\KK}) \in \SH(D)_{h \geq n}
\end{align*}
where the last connectivity statement is proved as in %
\cite{Spitzweck20}*{Lemma~5.5}, 
using $s_n \KGL_{\KK} = \Si^{2n,n} \MM\Z_{\KK}$. By \cite{Hoyois15}*{Lemma~2.2}, we have $i^* j_* (f_n \KGL_{\KK}) \in \SH(\kk)_{\geq n} = \SH(\kk)_{h\geq n}$.

For any $E \in \SH(\kk)_{\geq n}$ and $q \in \Z$, we have $f_q E \in \SH(\kk)_{\geq n}$, by \cite{Bachmann17}*{Proposition~4(3)}, and also $s_q E \in \SH(\kk)_{\geq n}$, because of the exact triangle $f_{q+1} E \to f_q E \to s_q E$. In particular, this yields
\[
s_{-1} i^* j_* (f_n \KGL_{\KK}) \in \SH(\kk)_{h\geq n}
\]
for all $n$. Since the connectivity of this tower goes to infinity as $n \to +\infty$, its homotopy limit is trivial %
\cite{Spitzweck20}*{Corollary~4.5}.
\end{proof}

\begin{lemma}\label{lem:ShriekNatural}
Let $E$ be a Cartesian motivic spectrum, $\xi \in \pi_{p,q} (\MM\F_{\ell} \sm E)$ a class which is compatible with base change, and $\tild{\xi} \colon \Si^{p,q} \MM\F_{\ell} \to \MM\F_{\ell} \sm E$ the adjunct map of $\MM\F_{\ell}$-modules. Then the purity isomorphism $i^! \MM\F_{\ell}^D \cong \Si^{-2,-1} \MM\F_{\ell}^{\kk}$ %
makes the following diagram in $\SH(\kk)$ commute:
\[
\xymatrix{
\Si^{-2,-1} \Si^{p,q} \MM\F_{\ell}^{\kk} \ar[d]^{\Si^{-2,-1} \tild{\xi}^{\kk}} & \Si^{-2,-1} i^* (\Si^{p,q} \MM\F_{\ell}^D) \ar[l]_{\cong} \ar[d]^{\Si^{-2,-1} i^* \tild{\xi}^D} \ar[r]^-{\pur_{\Si^{p,q} \MM}}_-{\cong} & i^! (\Si^{p,q} \MM\F_{\ell}^D) \ar[d]^{i^! \tild{\xi}^D} \\
\Si^{-2,-1} (\MM\F_{\ell}^{\kk} \sm E^{\kk}) \ar@{-}[d]^{\cong} & \Si^{-2,-1} i^* (\MM\F_{\ell}^D \sm E^{D}) \ar[l]_{\cong} \ar[r]^-{\pur_{\MM \sm E}} & i^! (\MM\F_{\ell}^D \sm E^D) \\
(\Si^{-2,-1} \MM\F_{\ell}^{\kk}) \sm E^{\kk} & (\Si^{-2,-1} i^* (\MM\F_{\ell}^D)) \sm i^* E^{D} \ar[l]_-{\cong} \ar[u]_-{\cong} \ar[r]^-{\pur_{\MM} \sm 1}_-{\cong} & i^! (\MM\F_{\ell}^D) \sm i^* E^D. \ar[u]_{\exmodr} \\
}
\]
\end{lemma}

Applying the lemma to the classes in the Milnor basis $\om(\seq) \in \pi_{*,*} (\MM\F_{\ell} \sm \MM\F_{\ell})$ yields:

\begin{corollary}\label{cor:ShriekSteenrod}
The purity isomorphism $i^! \MM\F_{\ell}^D \cong \Si^{-2,-1} \MM\F_{\ell}^{\kk}$ induces the commutative diagram in $\SH(\kk)$:
\[
\xymatrix{
\Si^{-2,-1} \bigoplus_{\seq} \Si^{p_{\seq, q_{\seq}}} \MM\F_{\ell}^{\kk} \ar[d]^{\Si^{-2,-1} \Psi^{\kk}} \ar[r]^-{\cong} & i^! (\bigoplus_{\seq} \Si^{p_{\seq, q_{\seq}}} \MM\F_{\ell}^{D}) \ar[d]^{i^! \Psi^D} \\
\Si^{-2,-1} (\MM\F_{\ell}^{\kk} \sm \MM\F_{\ell}^{\kk}) \ar[r] & i^! (\MM\F_{\ell}^{D} \sm \MM\F_{\ell}^{D}) \\
(\Si^{-2,-1} \MM\F_{\ell}^{\kk}) \sm \MM\F_{\ell}^{\kk} \ar@{-}[u]^{\cong} \ar[r]^-{\cong} & (i^! \MM\F_{\ell}^{D}) \sm \MM\F_{\ell}^{\kk}. \ar[u]_{\exmodr} \\
}
\]
\end{corollary}

\begin{remark}
In %
\cite{Spitzweck18}*{Theorem~10.26}, 
a \emph{right} $\MM\F_{\ell}$-module structure on $\MM\F_{\ell} \sm \MM\F_{\ell}$ is used. The arguments in Section~\ref{sec:ReduceLisse} also work in that case, with minor adjustments.
\end{remark}

\end{document}